\documentclass[11pt]{amsart}

\font\teneufm=eufm10
\font\seveneufm=eufm7
\font\fiveeufm=eufm5
\newfam\eufmfam
\textfont\eufmfam=\teneufm
\scriptfont\eufmfam=\seveneufm
\scriptscriptfont\eufmfam=\fiveeufm
\def\eufm#1{{\fam\eufmfam\relax#1}}
\def\Fgot{{\mathcal{F}} }
\def\Lgot{{\eufm L}}
\newcommand{\R}{\mathbb R}
\newcommand{\Q}{\mathbb Q}
\newcommand{\Z}{\mathbb Z}
\newcommand{\N}{\mathbb N}
\newcommand{\C}{\mathbb C}
\newcommand{\He}{\mathbb H}
\newcommand{\gra}{\nabla_{\He^n}}
\newcommand{\lap}{\Delta_{\He^n}}   
\newcommand{\BB}{{{\eufm{B}}}}
\renewcommand{\b}{\beta}
\newcommand{\s}{\sigma}
\newcommand{\z}{\zeta}
\newcommand{\e}{\varepsilon}
\newcommand{\cA}{{\mathcal{A}} }
\newcommand{\cV}{{\mathcal{V}} }

\newcommand{\cY}{{\mathcal{Y}} }
\newcommand{\om}{\omega}
\def\Tgot{{\mathcal{T}} }

\newtheorem{thm}{Theorem}[section]
\newtheorem{cor}[thm]{Corollary}
\newtheorem{lemma}[thm]{Lemma}
\newtheorem{prop}[thm]{Proposition}
\theoremstyle{definition}
\newtheorem{defn}[thm]{Definition}
\theoremstyle{remark}
\newtheorem{rem}[thm]{Remark}
\numberwithin{equation}{section}

\begin{document} 
\title[The Ginzburg-Landau equation
in
the Heisenberg group]{The Ginzburg-Landau equation
in the Heisenberg group}
\author{Isabeau Birindelli and Enrico Valdinoci}
\address{
Isabeau Birindelli,
Dipartimento di Matematica,
Universit\`a di Roma La Sapienza,
Piazzale Aldo Moro, 2,
I-00185 Roma (Italy)} 
\email{isabeau@mat.uniroma1.it}
\address{
Enrico Valdinoci,
Dipartimento di Matematica,
Universit\`a di Roma Tor Vergata,
Via della Ricerca Scientifica, 1,
I-00133 Roma (Italy)
}
\email{valdinoci@mat.uniroma2.it}%
\thanks{Partially supported by MURST Variational
Methods and Nonlinear Differential Equations.}

\begin{abstract}
We consider 
a functional related with 
phase transition models in
the Heisenberg group
framework. We prove that
level sets of local minimizers
satisfy some density estimates,
that is, they behave as ``codimension one" sets.
We thus deduce a uniform
convergence property of these level sets
to interfaces with minimal area.

These results are then applied in the construction
of (qua\-si)periodic,
plane-like minimizers, i.e., minimizers of our functional
whose level sets are contained in a spacial slab of
universal size in a prescribed direction.
As a limiting case, we obtain the existence of
hypersurfaces contained in such a slab which minimize
the surface area with respect to a given periodic metric.
\end{abstract}

\maketitle

\bigskip

\bigskip

\tableofcontents

\section*{Set up and notation}

Throughout this paper, 
$\chi_A$ will denote the charachteristic function 
of a set $A$, namely
$\chi_A(x)=1$ if $x\in A$ and $\chi_A(x)=0$
otherwise. Also,
the map $\R\ni \tau\mapsto \tau_+$ agrees with $\tau$
if $\tau\geq 0$ and with $0$ otherwise.

We will set $n\in\N$,
$n\geq 1$ and we consider the ``space''variables
$$ x=(x_1,\dots,x_n)\in \R^n\qquad{\mbox{and}}\qquad
y=(y_1,\dots,y_n)\in \R^n\,.$$
Also, 
we will often write $z=(x,y)\in\R^{2n}$
and we identify it with the vector $z\in\C^n$
with components $z_k=x_k+\sqrt{-1} \,y_k$, for $k=1,\dots, n$.
The vector $\overline z\in\C^n$
has components $\overline z_k=x_k-\sqrt{-1}\, y_k$,
for $k=1,\dots, n$.

Given $z,\,w\in\C^n$, the notation $zw$ stands
for the product in~$\C^n$, that is
$$ zw=\sum_{j=1}^n z_j w_j\,.$$

Given $t\in \R$,
we will use the notation
$\xi=(z,t)=(x,y,t)\in\R^{2n+1}$.
We will also consider the radial variables
$\rho'=|z|$ and 
\begin{equation}\label{radialvars} 
\rho= |\xi|_{\He^n}=
\Big( |z|^4+t^2\Big)^{1/4}\,.
\end{equation}
As usual, $\He^n$ denotes the Heisenberg group,
endowed with the action
$$ \xi\circ \xi_0\,=\,
\Big( 
z+z_0,\,t+t_0+2 \,{\rm Im}\,
(\overline z z_0)
\Big)\,.$$
By induction, one
sees that, if $K^{(1)},\dots,K^{(\ell)}\in\R^{2n}$,
then
\begin{eqnarray}\label{itekappa}
(K^{(\ell)},0)\circ\dots\circ (K^{(1)},0)\circ
(z,t) =\\
= \left(
z+\sum_{j=1}^\ell K^{(j)}\,,\;
t+2\,{\rm Im}\,\Big( \sum_{j=1}^\ell \overline K^{(j)} z
+\sum_{1\leq m<j\leq\ell} \overline K^{(j)} K^{(m)} \Big)
\right)\,
\nonumber.\end{eqnarray}
We will consider the so-called ``Koranyi
ball'' centered at $\xi_0$ of radius~$r$,
defined by
$$ \BB_r(\xi_0)\,=\,
\Big\{\xi \;\,
{\mbox{s.t.}}
\;\,|\xi^{-1}\circ \xi_0|_{\He^n}\leq r
\Big\}\,.$$
%% The notation $B_r(\xi_0)$ denotes
%% the standard Euclidean ball of radius~$r$
%% centered at~$x_0$.
We also set~$Q:=2(n+1)$.
If~$\Lgot^k$ denotes the $k$-dimensional
Lebesgue measure, one has that
$$ \Lgot^{2n+1} \Big( 
\BB_r(\xi_0)\Big)\,=\,
{\rm const} \,r^Q\,.$$
As usual in the
Heisenberg group setting, we consider the
vector fields which are left invariant with respect to the group action
$$ X_k =\partial_{x_k}+2 y_k\partial_t
\qquad{\mbox{and}}\qquad
Y_k =\partial_{y_k}-2 x_k\partial_t
\,,$$
for $k=1,\dots, n$ and which generate the whole Lie 
algebra together with 
their commutator. We also introduce the so-called
``Kohn Laplacian''
$$ \lap\,=\, \sum_{k=1}^n (X_k^2+Y_k^2)\,.$$
Given a (smooth) function $u$,
we also introduce the intrinsic gradient:
$$ \gra u\,=\, (X_1 u,\dots,X_n u,Y_1 u,\dots, Y_n u)\,.$$
For a radial (with respect to the Koranyi
ball)
function $v=v(\rho)$, 
a direct computation 
shows that
\begin{equation}\label{grast}
|\gra v|\,\leq\,{\rm const}\, |v'|
\end{equation}
and
\begin{equation}\label{lap:rad}
\lap v\,=\, \left(
\frac{|z|}{\rho}\right)^2\,
\left( v''+\frac{Q-1}{\rho}\,v'
\right)\,.
\end{equation}
See, for instance,
page~182 of
[B03] for further details.

\section{Statement of results}

Given a domain $\Omega\subseteq\R^{2n+1}$,
we define the functional
$$ \Fgot_\Omega(u)= 
\int_\Omega |\gra u(\xi)|^2+
F(\xi,u(\xi))\,d\xi\,,
$$
under the following structural assumptions
on the ``double-well potential''~$F$:
\begin{itemize}\label{7HHHGGHHhj11}
\item $F$ is non-negative and bounded and $F(\xi,1)=F(\xi,-1)=0$;
\item for any 
$\theta\in [0,1)$, $\displaystyle\inf_{|u|\leq\theta}
F(\xi,u)\geq\gamma (\theta)$, 
where $\gamma$
is decreasing and strictly positive;
\item there exist $d\in[0,2]$ and $\ell\in(0,1)$ so
that:
\begin{itemize}
\item $F(\xi,t)\geq{\rm const}\,(1-|t|)^d$, if $|t|\in(\ell,1)$;
\item $F$ is continuous in $u$, for $|u|<1$;
\item $F$ is locally Lipschitz continuous in~$u$
for~$|u|<1$;
\item if $d>0$, $F_u(\xi,u)$ is continuous for $|u|<1$, and,
if $s<\ell$, then
$$ F_u(\xi,-1+s)\geq {\rm const}\,s^{d-1}\,,\qquad
F_u(\xi,1-s)\leq -{\rm const}\,s^{d-1}\,;$$
\item if $d=2$, $F$ is continuous in $u$ for $|u|\leq 1$ and
$F_u(\xi,u)$ is 
increasing for $u$ near $\pm1$.
\end{itemize} \end{itemize}

%% Quantities depending only on $n$ and on the structural
%% contants here above will be referred to as
%% ``universal constants''.
Examples of such potentials\footnote{
We take this opportunity to amend a typo at the end of
page~164 of~[V04]: the condition
``$F_u(x,1-s)\geq -{\rm const}\,s^{d-1}$''
there
has to be
replaced 
by ``$F_u(x,1-s)\leq -{\rm const}\,s^{d-1}$''.}
are given by 
$F=Q(\xi)\,(1-u^2)^2$ and
$F=Q(\xi)\chi_{(-1,1)}(u)$.
The first type of potential,
when the Heisenberg group is replaced
by standard elliptic integrands,
is related to phase transition
and superfluid models
(see, e.g., [R79], 
[GL58] and~[GP58]).
The latter potential is related to flame propagation
and fluid jet models (see, e.g.,
[AC81] and~[ACF84]). 
The minimizers of~$\Fgot$ are lead to suitably approach
a step function, whose interface is a minimal surface 
(see~[G85] and~[M87] for details). 
Also, functionals of these 
type
have been the subject of a celebrated conjecture
of De Giorgi (see, for instance,~[DG79], [AAC01], [BL03]
and references therein).

We say that a function $u$ is a local minimizer for
$\Fgot$ in the
domain $\Omega$ if $\Fgot_\Omega (u)$ is well-defined and finite and
$$ \Fgot_\Omega(u+\varphi)\,\geq\,\Fgot_\Omega(u)\,,$$
for any $\varphi\in C^\infty_0 (\overline\Omega)$.

We prove the following density estimates, which,
roughly speaking, state that the
level sets of minimizers behave
``like sets of dimension $(Q-1)$''. 

\begin{thm}\label{CC2}
Fix $\delta>0$.
Let $u$ be a local minimizer for $\Fgot$ in a domain $\Omega$,
with $|u|\leq 1$.
Then:
\begin{itemize}
\item[(i)] there exist positive
constants $c$ and $r_0$, depending only on $\delta$
and on the structural constants, such that
$$ 
\Fgot_{\BB_r(\xi)}(u)\leq c\, r^{Q-1}
\,, $$
for any $r\geq r_0$,
provided that $\BB_{r+\delta}(\xi)\subseteq\Omega$;
\item[(ii)] for any $\theta_0\in(0,1)$,
for any $\theta \in[-\theta_0,\theta_0]$ and
for any $\mu_0>0$, if 
\begin{equation}\label{COND1}
\Lgot^{2n+1}
\Big(\BB_K(\xi)\cap\{u\geq \theta\}\Big)\,\geq\,
\mu_0\,,
\end{equation}
then there exist positive constants
$c^\star$ and $r_0$, depending
on $K$,
$\mu_0$, $\theta_0$ and on the structural constants,
such that
$$ \Lgot^{2n+1}
\Big( 
\BB_r(\xi)\cap \{u\geq \theta\}
\Big)\geq c^\star r^Q\,,$$
for any $r\geq r_0$,
provided that $\BB_{r+\delta}(x)\subseteq\Omega$.

Analogously, if
\begin{equation}\label{COND2}
\Lgot^{2n+1}
\Big(\BB_K(\xi)\cap\{u\leq \theta\}
\Big)\,\geq\,
\mu_0\,,\end{equation}
then
$$ \Lgot^{2n+1}
\Big(\BB_r(\xi)\cap \{u\leq \theta\}
\Big)\geq c^\star 
r^Q\,,$$
for any $r\geq r_0$,
provided that $\BB_{r+\delta}(\xi)\subseteq\Omega$.
\end{itemize}\end{thm}

The original idea of such density estimates
goes back to~[CC95] and several related techniques
have been recently developed in~[V04], [PV05a]
and~[PV05b]. 
In the Heisenberg
group setting,
a result analogous to the one in
Theorem~\ref{CC2}[{\rm i}]
was proven in~[BL03] for the case in 
which $F(\xi,u)=(1-u^2)^2$.
Also, our Theorem~\ref{CC2}
may be seen as the extension of analogous estimates
for minimal surfaces (see, e.g., [G84]).

As pointed out in [CC95], density estimates of
these type easily imply an $L^\infty_{\rm loc}$-convergence
of the level sets of the minimizers. Indeed, as proven in [MSC01],
if $u$ is a local minimizer and $u_\e(\xi)=u(z/\e,t/\e^2)$, then
$u_\e$ converges, up to subsequence, in $L^1_{\rm loc}$ to a step 
function, whith minimal interface (with respect to the $\He^n$-area). 
Then, using 
the argument in Theorem~2 in [CC95]
(see also the end of \S~6 of [PV04a]), one deduces from
Theorem~\ref{CC2}[{\rm ii}]
the
following result: 

\begin{thm}\label{unithm}
Fix $\theta\in(0,1)$.  Let $|u|\leq1$ be a local
minimizer of~$\Fgot$ in a bounded domain $\Omega$
and let
$u_\e(\xi)=u(z/\e,t/\e^2)$.  
Assume
that, as $\e$ tends to zero, $u_\e$ converges in
$L^1_{\rm loc}$ to
the step function $\chi_E - \chi_{\Omega\setminus E}$,
for a suitable $E\subset \Omega$. Then,
$\{|u_\e|\leq\theta\}$ converges locally uniformly to
$\partial E$, that is
$\mathrm{dist}\,(\xi,\partial E)\rightarrow 0$ uniformly for $\xi\in
\{|u_\e|\leq\theta\}\cap K$ for any set $K$ 
whose closure is contained in $\Omega$.
\end{thm} 

As a byproduct, we deduce from Theorem~\ref{unithm}
that the level sets of 
the 
example built
in [BL03], though different from being minimal surfaces themselves,
approach locally uniformly the hyperplane $\{t=0\}$, when
scaled via the natural $\He^n$-scaling.

Next is a version of Theorem~\ref{CC2}
which turns out to be convenient for applications:

\begin{thm}\label{IMPROVE}
Fix $\theta_0\in (0,1)$ and $\delta>0$.
Let~$u$ be a local minimizer
for~$\Fgot$ in a domain~$\Omega$, with $|u|\leq1$. Suppose that
$|u(\xi_0)|\leq \theta_0$.
Then, there exist positive constants~$c$ and~$r_0$,
depending only on~$\theta_0$, $\delta$
and on the structural constants, in such a way that
$$ \Lgot^{2n+1}\Big( \BB_r(\xi_0) \cap
\{|u|\leq \theta_0\}\Big)\,\geq\, c r^{Q-1}\,,$$
for any $r\geq r_0$, as long as $\BB_{r+\delta}(\xi_0)\subseteq
\Omega$.
\end{thm}

Another topic we deal with in this paper is related to
the construction of ``(qua\-si)periodic,
plane-like'' minimizers, that
is, local minima of~$\Fgot$ in any 
bounded domain
which are either periodic or locally uniformly
close to periodic, and whose level
sets are contained in a neighborhood
of a given hyperplane of universal size. 
These results extend the ones of~[V04]
to the Heisenberg
group and fit into the theory of plane-like
structures for PDEs in periodic media. They also generalize
the construction of orbits and minimal measures of prescribed rotation 
number (which is
a classic topic in
Aubry--Mather 
theory) to pseudoelliptic PDEs. We refer to~[M86], [B90], [CL98],
[CL01],
[RS03], [RS04], [CL05],
[PV05a] and [PV05b] for related results and further motivations.
We also recall that the study of the (quasi)periodic
level sets of suitable functions
is also a subject of enduring interest
in topology and mathematical physics 
(see, e.g.,~[DN05]).

We now present two results, namely 
Theorems~\ref{classa-per}
and~\ref{classa-qper}, according to whether we consider rational
or irrational frequencies.
More precisely, in the rational case,
we have
the following result:

\begin{thm}\label{classa-per}
Fix~$\delta>0$. Suppose that~$\om\in\Q^{2n}\setminus\{0\}$. 
Assume also
that
\begin{equation}\label{classperiodic}
F((k,0)\circ \xi,\cdot)=F(\xi,\cdot)
\end{equation}
for any $\xi\in\He^n$ and any $k\in\Z^{2n}$.
Then, there exists~$M_0>0$, depending only
on~$\delta$ and on the structural constants of~$F$,
such that the following holds.

There exists
a function $u_\om:\He^n\rightarrow [-1,1]$ 
which is a local minimizer for~$\Fgot$
in any bounded
domain of~$\He^n$ and so that
\begin{equation}\label{7uhnHHH}
\{ \xi\in \He^n \,,\;\,|u_\om|\leq 1-\delta \}\,\subseteq\,
\{ \xi\in \He^n \,,\;\,|\xi\cdot (\om,0)|\leq 
M_0\,|\om| \}\,.
\end{equation}
Moreover, $u_\om$ enjoys the following periodicity properties:
for any $\xi=(z,t)\in \He^n$ and any $k\in\Z^{2n}$ such that
$\om \cdot k=0$
\begin{equation}\label{7uhnHHHB}u_\om((k,0)\circ\xi)=u_\om
(\xi)\,,
\end{equation}
and
\begin{equation} \label{1.5bis}
u_\om(z,t)=u_\om(z, t+2)\,.
\end{equation}

Finally, $u_\om$ enjoys the following discrete
monotonicity property:
for any $k\in \Z^{2n}$ so that
\begin{equation}\label{7uhnHHHBi}{\mbox{
$\om\cdot k\geq 0$, we have that
$u((k,0)\circ \xi)\geq u(\xi)$ for any $\xi\in \He^n$.
}}\end{equation}
\end{thm}

\begin{rem} We would like to observe that 
the result of Theorem~\ref{classa-per}
cannot hold if
one fixes~$\om\in \Q^{2n+1}$ with $\om_{2n+1}\neq 0$.
Indeed, as
we will see
later, 
any function satisfying a periodicity condition as in~\eqref{7uhnHHHB}
turns out to be also 
periodic 
in the vertical variable and hence 
its level sets cannot be contained in a slab orthogonal to~$\om$
unless this slab is horizontal, that is, $\om_{2n+1}=0$.
This may be also rephrased by saying that
a function
cannot be monotone in a direction where it
is periodic. 
\end{rem} 

In the irrational frequency case, an analogous
statement holds up to locally uniform approximation. More precisely, the
following result holds true:

\begin{thm}\label{classa-qper}
Fix~$\delta>0$. Suppose 
that~$\om\in\R^{2n}\setminus\Q^{2n}$.
Assume
that~\eqref{classperiodic}
holds.
Then,
given any 
sequence of vectors $\om_j\in \Q^{2n}$ so that
$$\lim_{j\rightarrow +\infty } \om_j \,=\,\om\,,$$
there exists
a sequence of functions~$u_j:
\He^n\rightarrow[-1,1]$, which are
local minimizers for~$\Fgot$
in any bounded domain of~$\He^n$, which satisfy the level
set constraint and
the periodicity and monotonicity
properties in~\eqref{7uhnHHH},
\eqref{7uhnHHHB}, \eqref{1.5bis}
and~\eqref{7uhnHHHBi}
(with $\om_j$ replacing $\om$ there),
and which
converge to~$u_\om$ 
uniformly on compact subsets of~$\He^n$,
up to subsequences.
\end{thm}

We derive from 
Theorems~\ref{classa-per}--\ref{classa-qper}
a result on (quasi)periodic minimal surfaces
in the Heisenberg group, which
can be seen as the natural extension of
analogous results
poven in the Euclidean case by [M86] and [H32]
(for geodesics) and by [CL01] (for general
codimension~one minimal surfaces).

For this, we consider a bounded function
$\alpha:\He^n\rightarrow\R$. We suppose that
$$1/\alpha_0\geq \alpha(\xi)\geq\alpha_0>0$$ 
and that
$$\alpha((k,0)\circ \xi)=\alpha(\xi)$$
for any $\xi\in\He^n$ and $k\in\Z^{2n}$.
We consider the surface measure
induced by the metric~$\alpha$,
that is the surface measure in~$\He^n$
induced\footnote{More general metrics and wider class
of vector fields may alo be taken into account,
by performing computations similar to the ones we present
here.} 
by the vector field $\sqrt\alpha\gra$
(see [MSC01] and [FSSC01]
for details). Then,
given any~$\om\in \R^{2n}\setminus\{0\}$, we construct
a minimal surface at a universal distance
from the hyperplane $\{(\om,0)\cdot \xi=0\}$
which is either periodic or close to periodic surfaces,
according to whether $\om$ is rational or irrational.

More precisely, we prove the following result:

\begin{thm}\label{backto}
For any $\om\in\R^{2n}\setminus\{0\}$, 
there exists
a set $E_\om\subset\He^n$ in such a way 
$\partial E_\om$ is a local minimizer for
the surface measure in~$\He^n$
induced by the metric~$\alpha$ and
$\partial E_\om$ 
is contained in the slab
$$\{ \xi\in \He^n \,,\;\,|\xi\cdot (\om,0)|\leq 
M_0\,|\om| \}\,,$$
for a suitable~$M_0>0$, depending only 
on~$\alpha_0$ and on~$n$, but
independent of~$\om$.
\end{thm}

\section{Proof of Theorem \ref{CC2}[{\rm i}]}\label{sec:1}

We start with a variation of the standard Caccioppoli
inequality (as exposed, for instance, in [G03]).
Namely, let $\BB$ be a ball (say, centered at $0$)
of radius~$2$ and $\widetilde\BB$ be a concentric ball
of radius~$3$, and assume that both $\BB$ and $\widetilde\BB$
are well contained in $\Omega$.
Let $\tilde\eta\in C^\infty_0([0,3])$
be so that $0\leq \tilde\eta\leq 1$, $\tilde\eta(\tau)=1$
if $\tau\in [0,2]$ and $|\tilde \eta'|\leq 10$.
Let also $\eta(\xi)= \tilde\eta(\rho)$.
We observe that, by~(\ref{grast}),
$$ |\gra \eta|\,\leq\,{\rm const}\,.$$
Then, the minimality of~$u$ implies that
$$ \Fgot_{\widetilde\BB} (u)\,\leq\, 
\Fgot_{\widetilde\BB} (u-\eta^2 u)\,,$$
and therefore
\begin{eqnarray*}
0 &\leq& \int_{\widetilde\BB} \Big(
|\gra (\eta^2 u)|^2 -2\gra u\cdot\gra (\eta^2 u)
\Big) \,d\xi +\,{\rm const}\,\leq\\
 &\leq& \int_{\widetilde\BB} \Big(
4\eta ^2 u^2|\gra \eta|^2 
+\eta^4|\gra u|^2 + 4\eta^3 u\gra \eta\cdot\gra u-\\ &&
-2\eta^2|\gra u|^2
-4\eta u\gra \eta \cdot\gra u
\Big) \,d\xi +\,{\rm const}\,\leq\\
&\leq&-\int_{\widetilde\BB} 
\eta^2 \left(2-\eta^2 -\frac 1 2\right)|\gra u|^2\,d\xi
+\,{\rm const}\,,
\end{eqnarray*}
where we used a suitably scaled Cauchy Inequality
in the last estimate.
Thus, we obtain the following Caccioppoli-type
estimate:
\begin{equation}\label{cacci}
\frac 1 2\, \int_{\BB} |\gra u|^2
\leq\int_{\widetilde\BB} \frac{\eta^2}{2}|\gra u|^2\leq\,{\rm const}
\,.
\end{equation}
The proof of
Theorem~\ref{CC2}[{\rm i}]
is now a variation of the one
in~[CC95] (see also Proposition~3.2 of~[BL03]
and~[AAC01] for related techniques).
For short, in the sequel, we will omit in the notation
the centers of the balls considered.

Let $\tilde g\in C^\infty([0,\,r])$ be so that
$|\tilde g'|\leq 10$, 
$\tilde g(\tau)=-1$ if $\tau\in[0,\,r-1]$ and
$\tilde g(r)=1$ on $\partial B_r$.
We also introduce the radial (with respect to the
Koranyi ball) function
$g(\xi)=\tilde g(\rho)$.
By~(\ref{grast}),
$|\gra g|\leq\,{\rm const}$.

Define $u^\star=\min\{u,g\}$. Since $u$ is a local minimizer and
$$ \Lgot^{2n+1} (\BB_r\setminus\BB_{r-1})\leq\,{\rm const}\, r^{Q-1}
\,,$$ we have that
\begin{eqnarray}\label{est-i-}
\Fgot_{\BB_r}(u)&\leq& \Fgot_{\BB_r}(u^\star)\leq{\rm const}\,\left(
\int_{\BB_r\setminus \BB_{r-1}} |\gra u^\star|^2+ 
r^{Q-1}\right)\leq\nonumber
\\
&\leq&{\rm const}\,\left(
\int_{\BB_r\setminus
\BB_{r-1}} (|\gra u|^2+|\gra g|^2)+r^{Q-1}\right)\leq
\\ &\leq&{\rm const\,}\left(
\int_{\BB_r\setminus\BB_{r-1}} |\gra u|^2+r^{Q-1}
\right)\nonumber\,.\end{eqnarray}
Let us now cover $\BB_r\setminus \BB_{r-1}$ with 
Koranyi balls
$\BB_1,\dots
\BB_K$ of radius~$2$; by measure theoretic considerations, one 
sees that it is possible to take
$K\leq {\rm const}\,r^{Q-1}$. Let $\widetilde{\BB_j}$ be the ball
concentric to $\BB_j$ of radius~$3$.
Then, from~(\ref{est-i-}) and~(\ref{cacci}), we infer that
\begin{eqnarray*}
\Fgot_{\BB_r}(u)&\leq&
{\rm const\,}\left(
\sum_{j=1}^K\int_{\BB_j} |\gra u|^2+r^{Q-1}
\right)
\leq\\&\leq&
{\rm const\,}( K+
r^{Q-1})
\,,\end{eqnarray*}
which implies the desired result.

\section{Proof of Theorem \ref{CC2}[{\rm ii}]}\label{sec:2}

We prove the first claim, the second one being
analogous. Furthermore,
we observe that, with no loss of
generality, we may assume
$\theta$ as
close to $-1$ as we wish. Indeed: assume\label{closeto-1}
the result to be true for $\theta^\star$ (say, close to  
$-1$), and let $\theta\in[-\theta_0,\theta_0]$, with
$\theta^\star\leq-\theta_0$. Then
$$\mu_0\leq\Lgot^{2n+1}\Big( \{u\geq \theta\}\cap \BB_K\Big)\leq
\Lgot^{2n+1}\Big( \{u\geq\theta^\star\}\cap \BB_K\Big)\,,$$
therefore, using the result for $\theta^\star$,
we deduce from Theorem~\ref{CC2}[{\rm ii}], that
\begin{eqnarray*}
&& \,{\rm const\,}
r^{Q}\,\leq\\
&\leq &\Lgot^{2n+1}
\Big( \{u\geq \theta^\star\}\cap \BB_r\Big)\leq
\\ &\leq& \Lgot^{2n+1}
\Big( \{u\geq\theta\}\cap \BB_r\Big)+
\Lgot^{2n+1}
\Big( \{\theta^\star\leq u<\theta\}\cap \BB_r\Big)\leq\\
&\leq&\Lgot^{2n+1}
\Big( \{u\geq\theta\}\cap \BB_r\Big)+\frac{1}{
\displaystyle\inf_{u\in[\theta^\star,\theta_0]} F}\int_{\BB_r} 
F(\xi,u)\,d\xi\leq\\
&\leq&
\Lgot^{2n+1}
\Big( \{u\geq\theta\}\cap \BB_r\Big)
+\,{\rm const\,} \Fgot_{\BB_r} (u)\leq\\
&\leq&\Lgot^{2n+1}
\Big( \{u\geq\theta\}\cap \BB_r\Big)+{\,\rm const}
\, r^{Q-1}\,,
\end{eqnarray*}
which gives
$$\Lgot^{2n+1}
\Big( \{u\geq\theta\}\cap \BB_r\Big)\geq {\rm\, const\,}
r^Q\,.$$
Thus, in the rest of the proof, we may and
do assume that $\theta$ is appropriately
close to~$-1$.\medskip

Also, for further use, we recall the following
recursive results:

\begin{lemma}\label{rec-lemma1}
Fix $\nu\in\N$. 
Let $v_k\geq0$ and $a_k\geq 0$ be two
nondecreasing sequences such that $v_1+a_1\geq c_0$,
$$  v_k^{(\nu-1)/\nu}\leq
C_0\left(v_{k+1}+a_{k+1}-v_k-a_k-c_1a_k\right)
\,,$$
for any $k\in\N$ and some positive constants $c_0$, $c_1$, $C_0$.
Then, there exists
$\kappa>0$, depending on $c_0$, $c_1$,
and $C_0$ such that
$$  v_k+a_k\geq \kappa \,k^{\nu}
$$  for any $k\in\N$.
\end{lemma}

\begin{lemma}\label{rec-lemma2}
Fix $\nu\in\N$.
Let $a_k\geq 0$ be a sequence such that $a_1\geq c_0$, $a_k\leq C_0
\,L^\nu k^{\nu-1}$,
$$ \Big( \sum_{1\leq j\leq k} a_j\Big)^{(\nu-1)/\nu} \,\leq\, C_0
\,\Big(
a_{k+1}+\sum_{1\leq j\leq k} e^{-L(k+1-j)}a_j \Big),$$
for any
$k\in\N$ and some positive constants $L$, $c_0$, and $C_0$. Then, if
$L$ is suitably large (in dependence of~$c_0$ and~$C_0$), 
there exists
$\kappa>0$, depending on $c_0$ and $C_0$, such that
$$ a_k\geq \kappa\, k^{\nu-1}
$$  for any $k\in\N$.
\end{lemma}

For the proof of the above lemmata, see~[CC95] or, for further
details, Lemmata~2.1 and~2.2 in [PV05b].
Here, we just mention that these iterations may
be seen as a ``discrete version'' of a differential
estimate linking volumes and areas in the Isoperimetric
Inequality (see again~[CC95] for very nice heuristics
on this).
With this, we define the quantities which
will play here the roles of the volume and of the area.
Namely, we define
\begin{eqnarray*}
\cV(r)&=& \Lgot^{2n+1} (\BB_r\cap \{u>\theta\})\,,\\
\cA(r)&=&\int_{\BB_r\cap\{ u\leq\theta\}}F(\xi,u)\,d\xi\,.
\end{eqnarray*}

The proof of
Theorem \ref{CC2}[{\rm ii}]
follows the ideas of~[CC95], as developed in~[V04],
[PV05a] and~[PV05b], and it is divided into three cases, according
to whether $d=0$, $d\in(0,2)$ or $d=2$.
In any of these cases,
$A$ will denote a suitably large
free parameter, to be conveniently chosen.

\subsection{The case $d=0$}
For any $\tau\in [0,r]$, let
$$\tilde h(\tau)=
2(\tau-r+1)_+^{2}-1$$
and 
define the
following radial (with respect to the
Koranyi ball) function
\begin{equation}\label{rad:09871}
h(\xi)=\tilde h(\rho)\,.
\end{equation}
Let also
\begin{equation}\label{rad:09872}
\s\,=\,\min \{u, \,h\}
\qquad{\mbox{and}}\qquad
\b\,=\,\min\{u-\s,\, 1+\theta\}\,.\end{equation}
Exploiting~(\ref{lap:rad}), one sees that 
\begin{equation}\label{boundlapAA}
|\lap h|\leq
\,{\rm const}\end{equation} 
in $\BB_r\setminus\BB_{r-1}$. Moreover,
by the Cauchy and Sobolev Inequalities (for the latter,
applied to $\b^2$,
see, e.g., %%% formula (2.2) in [GN96]),
[VSCC92], pages~3 and~59, and the references
on page~62),
and the minimality of~$u$, we deduce that
\begin{eqnarray*}\label{sobca}
&& \Big( \int_{\BB_r} \b^{\frac{2Q}{Q-1}}
\Big)^{\frac{Q-1}{Q}}\nonumber\leq
\\ &\leq& 
{\rm \,const\,}
\int_{{\BB_r}\cap\{u-\s\leq
1+\theta\}}|\b|\,|\gra \b|\nonumber
\leq\\
&\leq& {\rm \,const\,}
A\Big( \int_{{\BB_r}\cap\{
u>\s\}}
(|\gra u|^2-|\gra\s|^2
-\\ &&-2\gra (u-\s)\cdot\gra\s )
\Big)+\nonumber \\ &&+
\frac{{\rm \,const\,}
}{A}\int_{{\BB_r}\cap\{u-\s\leq
1+\theta\}}(u-\s)^2\leq \\ &\leq&
{\rm \,const\,}
A\Big( \int_{{\BB_r}\cap\{u>\s\}}
(|\gra u|^2-|\gra \s|^2        
+\\ &&+2 (u-\s)\,\lap\s )
\Big)+
\frac{
{\rm \,const\,}
}{A}\int_{{\BB_r}\cap\{u-\s\leq
1+\theta\}}(u-\s)^2\leq \\ &\leq&
{\rm \,const\,}A\Big[
\int_{{\BB_r}\cap\{
u>\s\}} 
(F(\xi,\s)-F(\xi,u))
+\nonumber\\ &&+2\int_{{\BB_r}} (u-\s)\,\lap\s\Big]+
\frac{
{\rm \,const\,}
}{A}\int_{
\BB_r\cap\{u-\s\leq 
1+\theta\}}(u-\s)^2
\,.\end{eqnarray*}
Note that, up to now, the condition $d=0$
has not been used yet (this will allow
us to use similar procedures also for $0<d\leq2$).

We now observe
that the left hand side of the above
inequality is larger or equal than
$$ 
\Big( \int_{\BB_r \cap\{ u>\theta\}} \b^{\frac{2Q}{Q-1}}
\Big)^{\frac{Q-1}{Q}}\,\geq\, (1+\theta)^2 
\,\cV(r-1)^{\frac{Q-1}{
Q}}\,.$$
By our assumptions on~$F$, the contribution
of the right hand side in~$\BB_{r-1}$ (where~$\s=-1$)
is bounded by
\begin{eqnarray*}&&
-\,{\rm const}\,A \int_{\BB_{r-1}} F(\xi,u) +\frac{\,{\rm const}\,}{A}
\int_{\BB_{r-1}\cap \{u\leq\theta\}}
(u+1)^2\leq\\ &\leq&
- \int_{\BB_{r-1}\cap \{u\leq\theta\}} \left[ 
{\,{\rm const}\,}{A}
-
\frac{\,{\rm const}\,}{A}
\right]\, F(\xi,u)\,d\xi\leq \\&\leq& -
{\,{\rm const}\,}
\cA (r-1)\,,
\end{eqnarray*}
if $A$ is conveniently large.
Moreover, exploiting~(\ref{boundlapAA}) and the
fact that $d=0$, one bounds
the contribution of the right
hand side in $\BB_r\setminus\BB_{r-1}$
by \begin{eqnarray*}
&& {\,{\rm const}\,}\Big[
\int_{(\BB_r\setminus\BB_{r-1})\cap \{u\leq \theta\}}
\Big(F(\xi,\s)-F(\xi,u)+\\ &&+(u-\s)+(u-\s)^2\Big)\,d\xi +\cV(r)-\cV(r-1)
\Big]\leq\\
&\leq&  \,{\rm const}\,\left[
\int_{(\BB_r\setminus\BB_{r-1})\cap \{u\leq \theta\}}
F(\xi,u)\,d\xi+\cV(r)-\cV(r-1)
\right]
\leq\\ 
&\leq&
{\,{\rm const}\,}
\Big[
\cA(r)-\cA(r-1)+
\cV(r)-\cV(r-1)
\Big]
\,.
\end{eqnarray*}
Then,
\begin{eqnarray*} &&
{\rm \,const\,}\Big( (\cV(r-1))^{\frac{Q-1}{Q}} 
+\cA(r
-1)\Big)\leq\\ &\leq&
\cV(r)-\cV(r-1)+\cA(r)-\cA(r-1)\,,\end{eqnarray*}
which completes the proof of the desired result
via Lemma~\ref{rec-lemma1}.

\subsection{The case $d\in(0,2)$}

Let $d'=\max\{1, d\}$ and, for $\tau\in[r-1,r]$,
we define 
$$\tilde h(\tau)\,=\,
2(\tau-r+1)_+^{\frac{2}{2-d'}}-1\,.$$
Let $h$, $\b$ and $\s$ be as in~(\ref{rad:09871})
and~(\ref{rad:09872}).
Recalling~(\ref{lap:rad}),
%%% and the fact that $d+2\geq 2d'$, 
one obtains 
that, in $\BB_r\setminus\BB_{r-1}$,
$$ |\lap h|\leq {\,\rm const\,}
(h+1)^{d'-1}\,.$$
Also, in the set $\{u>\s\}$ we have $\s=h$ and then
the last formula implies
\begin{eqnarray}\label{jhschss}
\lap h\, (u-\s)&\leq&
{\,\rm const\,}(h+1)^{d'-1}\,
(u-h)\leq\nonumber\\
&\leq& {\,\rm const\,}(u+1)^{d'}\,,\end{eqnarray}
in $(\BB_r\setminus\BB_{r-1})\cap \{u>\s\}$.

As argued in the case $d=0$,
Sobolev (applied to $\b^2$)
and Cauchy Inequalities
and the minimality property
of~$u$ imply that
\begin{eqnarray}\label{eq***}
&& \cV(r-1)^{\frac{Q-1}{Q}}\leq
\nonumber\\ &\leq&
{\,\rm const\,}
A\Big[
\int_{\BB_r}
(F(\xi,\s)-F(\xi,u))
+\nonumber\\ &&+2\int_{\BB_r}
\lap\s\,(u-\s)
\Big]+\\ &&+\nonumber\frac{
{\,\rm const\,}
}{A}\int_{
\BB_r\cap\{u-\s\leq
1+\theta\}}(u-\s)^2
\,.\end{eqnarray}
We now estimate the right hand side
in~$\BB_{r-1}$, in which~$\s=-1$.
For this, note that, if $\theta$ is close to~$-1$ and
$u\leq \theta$, we have, by our assumptions on~$F$, that
\begin{equation}\label{small}
F(\xi,u)\geq {\,\rm const\,}
(1+u)^d\geq\frac{
{\,\rm const\,}
}{2^{2-d}}(1+u)^2\,
\end{equation} so that
$$\frac{
{\,\rm const\,}
}{A}\int_{\BB_{r-1}\cap\{u\leq
\theta \}}(1+u)^2-\frac{
{\,\rm const\,}
A}{2}\int_{\BB_{r-1}\cap\{u\leq
\theta \}}
F(\xi,u)\leq 0\,$$
if we choose $A$ suitably big. Therefore,
the right hand side of (\ref{eq***})
in~$\BB_{r-1}$ is negative and less than
$-{\,\rm const\,}\,\cA(r-1)$.

Note also that,
since $\theta$ is close to $-1$ and $\s\leq u$,
our assumptions on
$F$ imply that $F(\xi,\s)\leq F(\xi,u)$ if $u\leq\theta$.
Thus, 
\begin{eqnarray*}&&
{\,\rm const\,}A \int_{\BB_r\setminus\BB_{r-1}}F(\xi,\s)-F(
\xi,u)\leq\\ &\leq&
{\,\rm const\,}A\,\Lgot^{2n+1}(\{ u>\theta\}\cap
(\BB_r\setminus\BB_{r-1}) )
\,.
\end{eqnarray*}
What is more, using again~(\ref{small}), one sees that
\begin{eqnarray*}
&&\frac{
{\,\rm const\,}
}{A}\int_{\BB_r\setminus\BB_{r-1}}(u-\s)^2\leq
\\ &\leq& \frac{
{\,\rm const\,}
}{A}\int_{(\BB_r\setminus\BB_{r-1})\cap
\{u\leq \theta\}
}(u+1)^2+\\ &&+\frac{
{\,\rm const\,}
}{A}\Lgot^{2n+1}(\{ u>\theta\}\cap (\BB_r\setminus
\BB_{r-1}))
\leq
\\ &\leq&
\frac{
{\,\rm const\,}
}{A}\int_{(\BB_r\setminus\BB_{r-1})\cap\{ u\leq\theta
\}}F(\xi,u)+\\
&&+\frac{
{\,\rm const\,}
}{A}\Lgot^{2n+1}(\{ u>\theta\}\cap (\BB_r\setminus\BB_{r-1}))\,.
\end{eqnarray*}
Also, recalling that $\s=-1$ in~$\BB_{r-1}$,
splitting the domain of integration into the sets
$\{u\leq\theta\}$ and $\{u>\theta\}$, using (\ref{jhschss})
and our hypotheses on~$F$,
we get
\begin{eqnarray*}&&
\int_{{\BB_r}}
%% \Big(
\lap\s\,(u-\s)
%% +\gra\s\,(u-\s)\Big) 
=\\
&& =
\int_{{(\BB_r\setminus\BB_{r-1})\cap \{u>\s\} }}
%% \Big(
\lap\s\,(u-\s)
%% +\gra \s\,(u-\s)\Big)
\leq\\
&&\leq {\,\rm const\,}\Lgot^{2n+1}\Big( 
(\BB_r\setminus\BB_{r-1})\cap\{u>\theta\}\Big)
+\\&&+
{\,\rm const\,}\int_{(\BB_r\setminus\BB_{r-1})\cap\{u\leq\theta\}} 
(u+1)^{d'}\leq\\
&&\leq {\,\rm const\,}\Lgot^{2n+1}\Big( 
(\BB_r\setminus\BB_{r-1})\cap\{u>\theta\}\Big)
+\\&&
+{\,\rm const\,}\int_{(\BB_r\setminus\BB_{r-1})\cap\{ u\leq\theta
\}} F(\xi,u)\,.
\end{eqnarray*}
Collecting the above estimates and using Lemma~\ref{rec-lemma1},
the result follows.

\subsection{The case $d=2$}

We will use here two further positive 
free parameters $\Theta$ and $T$: we will fix 
$\Theta$ small enough and then choose $T$ so that
$\Theta T$ is suitably large.

Let $k\in\N$. On page~183 of [V04],
a function $\tilde h\in C^2 ([0,(k+1)T])$
was constructed so that $-1\leq \tilde h\leq 1$,
$\tilde h((k+1)T)=1$, $\tilde h'(0)=0$,
\begin{equation}\label{veca}
\tilde h(\tau)+1\,\leq\,{\rm const}\, e^{-\Theta T(k+1-j)}
\end{equation}
if $\tau\in [(j-1)T,\,jT]$, for $j=1,\dots,k+1$,
$$ |\tilde h'(\tau)|\leq {\,\rm const\,}\Theta \tau (\tilde h+1)
$$
if $\tau\in [0,1]$,
$$ |\tilde h'(\tau)|\leq {\,\rm const\,}\Theta  (\tilde h+1)
$$
if $\tau\in [1,(k+1)T]$, and
$$
|\tilde h''(\tau)|\leq {\,\rm const\,}\Theta (\tilde h+1)
$$
if $\tau\in [0,(k+1)T]$. 

Let $h$ be as in~(\ref{rad:09871}).
In the light of~(\ref{lap:rad}), we thus have that
\begin{equation}\label{-1-}
|\lap h|\leq 
{\,\rm const\,}
(h+1)\leq {\,\rm const\,}
e^{-\Theta T(k+1-j)}\end{equation}
in $\BB_{jT}\setminus\BB_{(j-1)T}$, and
\begin{equation}\label{-1'-}
|\lap h|\leq {\,\rm const\,}\Theta(h+1)
\end{equation}
in $\BB_{(k+1)T}$. Consequently, 
if~$\Theta$ is suitably small
and~$h$ is suitably close to~$-1$, then
\begin{equation}\label{-2-}
|\lap h|\leq \sqrt{\Theta} F_u (x,h)
\end{equation}
in $\BB_{(k+1)T}$. 
We set \begin{equation}
\label{defcsta}\theta'=\theta-C_*
e^{-\Theta T}\,,\end{equation}
where $C_*$ denotes a suitably large constant.
Note that, if $\Theta T$ is large 
enough, then
$\theta'>-1$. 
Define also $$
\s\,=\,\min\{u,h\}\quad{\rm and}\quad
\b\,=\,\min\{u-\s,\,1+\theta'\}\,.$$
By arguing as done on page~\pageref{sobca},
we
gather from the Sobolev and Cauchy Inequalities that
\begin{eqnarray}\label{-3-}
&&\Big( \int_{\BB_{(k+1)T}}
\b^{\frac{2Q}{Q-1}}
\Big)^{\frac{Q-1}{Q}}\leq \nonumber\\
&\leq&{\,\rm const\,}
A\left[ \int_{\BB_{(k+1)T}\cap \{ u>\s
\} }\Big(
F(\xi,\s)-F(\xi,u)+\lap\s(u-\s)\Big)\right]+\nonumber\\ 
&&+\frac{{\,\rm const\,}}{A}
\int_{{\BB_{(k+1)T}\cap \{ u-\s\leq
1+\theta' \} }}
(u-\s)^2\,.
\end{eqnarray}
We now estimate the left hand side
of (\ref{-3-}). If $\Theta T$ is large enough,
we gather from~(\ref{veca}) that
$\theta-h\geq (1-\theta_0)/2$ in $\BB_{kT}$. 
Thus,
the left hand side
of (\ref{-3-}) is bigger than 
$${\,\rm const\,}\cV(kT)^{\frac{Q-1}{Q}}\,.$$
Let us now estimate the right hand side of (\ref{-3-}).
First of all, we consider the contribution in
$\{u\leq\theta\}$.
Let us observe that, since $-1\leq\s\leq u\leq 1$,
\begin{eqnarray*}
&& (u+1)^2-(\s+1)^2 -\frac{1}{2}(u-\s)^2=\\ 
%%% &=& (u+\s)\,(u-\s)+2(u-\s)-\frac{1}{2}(u-\s)^2   =\\ 
&=&(u-\s)\,\left(\frac{1}{2} u+\frac32\s+2\right)\geq 0\,.
\end{eqnarray*}
Accordingly,
in $\{\s<u\leq\theta\}$,
\begin{eqnarray}\label{-pound-}
&& F(\xi,u)-F(\xi,\s)=\int_h^u F_u(\xi,\z)\,d\z\geq\nonumber\\
&& \geq {\,\rm const\,}
\int_h^u(\z+1)\,d\z={\,\rm const\,}\Big[ (u+1)^2-(h+1)^2\Big]
\geq\nonumber\\&& \geq {\,\rm const\,}(u-h)^2\,.
\end{eqnarray}
Consequently, 
recalling again~(\ref{-2-}),
we deduce that
the contribution of the right hand side of (\ref{-3-})
in $\{u\leq\theta\}$ is controlled by
\begin{equation}\label{-4-}
\int_{\BB_{(k+1)T}\cap \{ \s<u\leq
\theta\}}\Big(
F(\xi,\s)-F(\xi,u)+{\,\rm const\,}\sqrt{\Theta} F_u(
\xi,\s)(u-\s)\Big)\,,
\end{equation}
as long as~$A$ is sufficiently large.

We now show that this quantity is indeed negative.
Since we assumed $\theta$ to be close to $-1$
(recall the discussion on page~\pageref{closeto-1}),
we have that $F$ and $F_u$ are monotone in 
%%%%% $(-1,\theta)$ and so, in 
$\{ \s<u\leq
\theta\}$, we have that $F(\xi,\s)-F(\xi,u)$ is negative
and that
$$ |F_u (\xi,\s)(u-\s)|\leq |F(\xi,\s)-F(\xi,u)|\,.$$
Since we assumed $\Theta$ to be small, 
the latter inequality yields
that
the quantity in~(\ref{-4-}) is negative.

Let us now bound
the right hand side of (\ref{-3-})
in $\{u>\theta\}$. 
The contribution in $\BB_{(k+1)T} \setminus\BB_{kT}$
of such term is bounded 
by
$$ \int_{(\BB_{(k+1)T} \setminus\BB_{kT})\cap \{u>\theta\} }
\Big(
F(\xi,\s)-F(\xi,u)+(\s+1) (u-\s) +(u-\s)^2
\Big)\,,$$
due to (\ref{-1'-}),
and this term is
bounded by
\begin{eqnarray*}&&\Lgot^{2n+1} \Big( \{u>\theta\}\cap 
(\BB_{(k+1)T}\setminus\BB_{kT})
\Big)\\
&=& \cV( (k+1)T )-\cV(kT)
\,.\end{eqnarray*}
Let us now look at the contribution
of the right hand side of
(\ref{-3-}) in $\{u>\theta\} \cap \BB_{k T}$. 
Notice that, from (\ref{veca}),
$$ \BB_{kT}\cap \{ \s<u\leq
\s+1+\theta'\}\subseteq
\BB_{kT}\cap \{ \s<u\leq
\theta\}\,,$$
provided that $C_*$ in (\ref{defcsta}) is large enough.
Consequently, the contribution
of the right hand side of (\ref{-3-}) in 
$\{u>\theta\} \cap \BB_{k T}$ is controlled by
\begin{eqnarray}\label{lalla}
& &\int_{\BB_{kT}\cap \{u>\theta\}}F(\xi,\s)-F(\xi,u)
+|\lap h|\leq\nonumber\\
&\leq& \sum_{j=1}^k \int_{\BB_{jT}\setminus\BB_{(j-1)T}\cap
\{u>\theta\}}
F(\xi,h)+|\lap h|\,.
\end{eqnarray}
By our assumption on $F$ in the case $d=2$, 
we have that
$$F(\xi,-1+s)\leq 
\,{\rm const}\,
s\,,$$
provided that $s>0$ is small enough. 
Thus, by means of~(\ref{veca}),
and~(\ref{-1-}),
we bound the above term in (\ref{lalla}) by
$$ \sum_{j=1}^k e^{{-\Theta T  (k+1-j)}}
\Big[\cV(jT)-\cV((j-1)T)\Big]\,.$$
Collecting all theses estimates, we get that
\begin{eqnarray*} 
&& {\rm\, const}\,(\cV(kT))^{\frac{Q-1}{Q}} 
\leq\\ &\leq&
\cV((k+1) T)-\cV(kT)\,+\,
\sum_{j=1}^k e^{{-\Theta T  (k+1-j)}}
\Big[\cV(jT)-\cV((j-1)T)\Big]
\,.\end{eqnarray*}
Then, the desired result follows from Lemma~\ref{rec-lemma2}.

\section{Proof of Theorem~\ref{unithm}}\label{sec:3}

This proof is a variation of the one on page~11 of~[CC95]
(see also page~69 of~[PV05a]). Due to some
subtleties given by the Heisenberg group
and the role
played by the Koranyi
ball, we provide full details for the facility of the reader.

The proof is based on Theorem~\ref{CC2}.
In fact, we first notice that conditions~(\ref{COND1})
and~(\ref{COND2}) in Theorem~\ref{CC2}, which are
somewhat unpleasant to check, may be
replaced by pointwise conditions (that is,
conditions~(\ref{cond=pt1})
and~(\ref{cond=pt2})
here below), that
are easier to deal with. Namely,
we deduce from Theorem~\ref{CC2} the following
result:

\begin{cor}\label{point-cor} 
Let $\delta$, $\Omega$ and $u$ be as in the
statement of Theorem~\ref{CC2}. Fix $\theta_0\in (0,1)$ and let
$\theta_1,\;\theta_2
\in (-\theta_0,\theta_0)$. Let $\xi\in \R^{2n+1}$ be so that
\begin{equation}\label{cond=pt1}
u(\xi)\geq \theta_1\,.
\end{equation}
Then, there exists positive constants $c^\star$ and $r_0$,
possibly depending on $\theta_0$ and other 
structural constants, such that
$$ \Lgot^{2n+1}\Big( \BB_r(\xi)\cap 
\{u\geq\theta_2\}\Big)\geq c^\star r^Q\,,$$
for any $r\geq r_0$, provided that $\BB_{r+\delta}(\xi)
\subseteq \Omega$.
Analogously, if
\begin{equation}\label{cond=pt2}
u(\xi)\leq \theta_1\,.
\end{equation}
then
$$ \Lgot^{2n+1}\Big( \BB_r(\xi)\cap 
\{u\leq\theta_2\}\Big)\geq c^\star r^Q
\,.$$
\end{cor}

\begin{proof}
We prove the first claim, the
second one being analogous.
We recall that $u$ is uniformly H\"older continuous
(see Theorem~2.1 in~[M95]), 
that is, there exist positive universal 
constants
$\Lambda$ and $\alpha$ so that
\begin{equation}\label{HhhHuUU} 
|u(\zeta)-u(\zeta')|\,\leq \,\Lambda
\,|(\zeta')^{-1}\circ\zeta|^\alpha\,,
\end{equation}
for any $\zeta,\;\zeta'\in \Omega'$, for $\Omega'$
well contained in $\Omega$.
We set
$$\hat\theta_1\,=\,
\frac{\theta_1-1}{2}\in \left[
\frac{-\theta_0-1}{2},
\frac{\theta_0-1}{2}\right]\subset(-1,1)
$$
and
$$ K\,=\,\left( 
\frac{1-\theta_1}{2\Lambda
}\right)^{1/\alpha}\,.$$
Thence, we conclude 
from~(\ref{cond=pt1})
and~(\ref{HhhHuUU}) that
$u\geq\hat\theta_1$ in $\BB_K(\xi)$.
Consequently,
$$
\Lgot^{2n+1}\Big( \BB_K(\xi)\cap 
\{u\geq\hat\theta_1\}\Big)=
\Lgot^{2n+1}\Big( \BB_K(\xi)\Big)=
{\,\rm const\,}K^Q\,,$$
thus condition~(\ref{COND1}) is fulfilled.
So, by Theorem~\ref{CC2}[{\rm ii}],
\begin{equation} 
\label{kkKKkUUu}
\Lgot^{2n+1}\Big( \BB_r(\xi)\cap 
\{u\geq\hat\theta_1\}\Big)
\geq {\,\rm const\,} r^Q\,.\end{equation}
Now: if $\hat\theta_1\geq \theta_2$,
the desired result follows 
directly from~(\ref{kkKKkUUu}).
If, on the other hand, $\hat\theta_1<\theta_2$,
we combine~(\ref{kkKKkUUu}) with
Theorem~\ref{CC2}[{\rm i}], to conclude that
\begin{eqnarray*}
&&
\Lgot^{2n+1}\Big( \BB_r(\xi)\cap 
\{u\geq\theta_2\}\Big)\geq\\ &\geq&
\Lgot^{2n+1}\Big( \BB_r(\xi)\cap 
\{u\geq\hat\theta_1\}\Big)-
\Lgot^{2n+1}\Big( \BB_r(\xi)\cap 
\{\hat\theta_1 \leq u\leq\theta_2\}\Big)
\geq\\ &\geq&
\Lgot^{2n+1}\Big( \BB_r(\xi)\cap 
\{u\geq\hat\theta_1\}\Big)-
{\,\rm const\,} \Fgot_{\BB_r(\xi)}(u)\geq\\
&\geq& 
{\,\rm const\,} r^Q -{\,\rm const\,}r^{Q-1}\geq\\
&\geq& {\,\rm const\,} r^Q\,,
\end{eqnarray*}
if $r$ is large enough.
\end{proof}

With this result, we may now complete the
proof of Theorem~\ref{unithm}, by arguing
as follows.

The proof 
of Theorem~\ref{unithm}
is by contradiction. If the claim of Theorem~\ref{unithm}
were false, there would exist a compact set $K$ (say, contained
in the Korany ball of radius~$R$), well contained in $\Omega$,
an infinitesimal sequence of~$\e_k$ and a positive~$\delta$ so that
$$ \sup_{\xi \in \{ |u_{\e_k}|\leq \theta\} \cap K}
\mathrm{dist}\,(\xi,\partial E)\,\geq\,8\delta\,,$$
for any $k \in\N$.
Consequently, we find a sequence of points
$\xi_k \in \{ |u_{\e_k}|\leq \theta\} \cap K$, in such
a way that
\begin{equation}\label{dist-EQ}
\mathrm{dist}\,(\xi_k,\partial E)\,\geq\,4\delta^2\,.
\end{equation}
Since, by construction, $\xi_k\in K\subseteq 
\BB_R(0)$,
we also have that
\begin{equation}\label{A-00}
|\xi_k|_{\He^{n}}\,\leq\,R\,.
\end{equation}
Of course, by possibly reducing $\delta$,
we may also assume that
\begin{equation}\label{B-00} \delta \,<\,\min\left\{ 
1\,,\;\frac{ \mathrm{dist}\,(K,\partial 
\Omega)}{16}\right\}\,,
\end{equation}
therefore~(\ref{dist-EQ}) implies that either
\begin{equation}\label{poss1}
\BB_{2\delta}(\xi_k) \subseteq \Omega\cap E 
\end{equation}
or
\begin{equation}\label{poss2}
\BB_{2\delta}(\xi_k) \subseteq \Omega\setminus E 
\end{equation}
for infinitely many $k\in \N$.
We assume that~(\ref{poss1}) holds (the case in which~(\ref{poss2})
holds is indeed analogous). Then,~(\ref{poss1})
implies that
\begin{equation}\label{chi=1}
\chi_E - \chi_{\Omega\setminus E}\,=\,1\qquad\quad{\mbox{in 
$\BB_{2\delta}(\xi_k)$}}
\end{equation}
for infinitely many $k\in \N$.
Also, by construction, if we set $\xi_k=(z_k,t_k)$ and
$\z_k= (z_k/\e_k, t_k/\e_k^2)$ we have that
$|u(\z_k)|\leq \theta$. Therefore, 
%%% if we define
%%% $$\theta':= \frac{1+\theta}{2}\,\in\, (\theta,1)
%%% \,,$$
we conclude from Corollary~\ref{point-cor}
that
$$ \Lgot^{2n+1} 
\Big(
\BB_r(\z_k)\cap \{u\leq \theta\}
\Big)\geq \,{\rm const}\,r^Q\,,$$
for any $r\geq r_0$, with $r_0$ universal.
Scaling back, this means that
\begin{equation}\label{ball=1} \Lgot^{2n+1}
\Big(
\BB_{\e_k r}(\xi_k)\cap \{u_{\e_k} \leq \theta\}
\Big)\geq \,{\rm const}\,(\e_k r)^Q\,,
\end{equation}
for any $r\geq r_0$.
We now set
$$ \delta'= \frac{\delta}{2(1+8R)}\,.$$
From~(\ref{chi=1}) and~(\ref{ball=1}), we conclude that
\begin{equation}\label{C-00}
\chi_E - \chi_{\Omega\setminus E}\,=\,1\qquad\quad{\mbox{in
$\BB_{\delta'}(\xi_k)$}}\end{equation}
and that, from (\ref{ball=1}) with $r=\frac{\delta'}{\e_k}$,
\begin{equation}\label{D-00}
\Lgot^{2n+1}
\Big(
\BB_{\delta'}(\xi_k)\cap \{u_{\e_k} \leq \theta\}
\Big)\geq \,{\rm const}\,(\delta')^Q\,,
\end{equation}
for infinitely many $k\in\N$.
Also,
if $$ \Omega' \,=\,
\bigcup_{k\in\N}\BB_{\delta'} (\xi_k)\,,$$
we deduce from~(\ref{B-00}) that $\Omega'$ is
well contained in $\Omega$.
We utilize this fact, together with~(\ref{C-00})
and~(\ref{D-00}), to infer that
\begin{eqnarray*}
&& \int_{\Omega'} |u_{\e_k} - (
\chi_E - \chi_{\Omega\setminus E})|\geq\\ &\geq&
\int_{
\BB_{\delta'}(\xi_k)\cap \{u_{\e_k}\leq \theta\}
} |u_\e -( \chi_E - \chi_{\Omega\setminus E})|\geq\\
&\geq&
\int_{
\BB_{\delta'}(\xi_k)\cap \{u_{\e_k} \leq \theta\}
} (1-\theta)\geq\\ &\geq& {\rm const}\,(\delta')^Q\,
({1-\theta})\,,
\end{eqnarray*}
for infinitely many $k\in\N$. This contradicts
the assumption that $u_{\e_k}$ converges to
$\chi_E - \chi_{\Omega\setminus E}$ in $L^1_{\rm loc}$,
thus proving Theorem~\ref{unithm}.

\section{Proof of Theorem~\ref{IMPROVE}}\label{sec:4}
Let $\theta_1:= (\theta_0+1)/2\in (\theta_0,1)$.
Then, since~$u$
is uniformly continuous, due to~[M95],
we get that
$$ \BB_{\e_0} (\xi_0)\,\subseteq\,
\{ |u|\leq\theta_1\}\,,$$
for a universal~$\e_0>0$.
Consequently,
$$ \Lgot^{2n+1} ( \BB_{\e_0} (\xi_0) \cap \{ u\geq -\theta_1\})
\geq \Lgot^{2n+1} ( \BB_{\e_0} (\xi_0) \cap \{ |u|\leq \theta_1\})
\geq \,{\rm const}\, \e_0^Q$$
and therefore, by Theorem~\ref{CC2}[{\rm ii}],
\begin{equation}\label{8iujn7yh34ed1}
\Lgot^{2n+1} ( \BB_{r} (\xi_0) \cap \{ u\geq -\theta_1\})
\,\geq\,{\rm const}\,r^Q\,,
\end{equation}
for~$r$
as in the statement of Theorem~\ref{IMPROVE}.
On the other hand, using Theorem~\ref{CC2}[{\rm i}],
we have that
\begin{equation}\label{8iujn7yh34ed2}
\Lgot^{2n+1}  ( \BB_{r} (\xi_0) \cap \{ -\theta_1 \leq u\leq
\theta_0\})\leq \,{\rm const}\,r^{Q-1}\,.
\end{equation}
By collecting the estimates in~\eqref{8iujn7yh34ed1}
and~\eqref{8iujn7yh34ed2}, we deduce that
\begin{equation}\label{23ed5551}
\begin{split}
\Lgot^{2n+1}  ( \BB_{r} (\xi_0) \cap \{ u>
\theta_0\})\,&\geq\, {\rm const}\,r^{Q}\,-
\,{\rm const}\,r^{Q-1}\,\geq
\\ \,&\geq\,\,{\rm const}\,r^{Q}\,,
\end{split}\end{equation}
as long as~$r_0$ is large enough.

Analogously,
\begin{equation}\label{23ed5552}
\Lgot^{2n+1}  ( \BB_{r} (\xi_0) \cap \{ u<
-\theta_0\})\,\geq\, {\rm const}\,r^{Q}\,.
\end{equation}

Let now~$S^{k,p}$
denote the Sobolev-type space in the Heisenberg
group setting (see, for instance,~[FS74]). Since~$u$
is in~$S^{1,2}(\BB_r)$ by construction
(and so in~$S^{1,1}(\BB_r)$)
and it is continuous by~[M95], we can find a
smooth
function~$v$ in such a way
\begin{equation}\label{8uijnhjHHH}
\int_{\BB_r(\xi_0)} |\gra u-\gra v|\,\leq\, 1
\end{equation}
and
\begin{equation}\label{8uijnhjHHH2}
{\mbox{
$|u(\xi)-v(\xi)|\leq \displaystyle\frac{\theta_0}4$ for any $\xi\in
\BB_r(\xi_0)$.}}\end{equation}
Let also
$$ w(\xi):=\left\{\begin{matrix}
v(\xi) & {\mbox{ if $|v(\xi)|\leq\theta_0/2$,}}\\
\theta_0/2 & {\mbox{ if $v(\xi)\geq\theta_0/2$,}}\\
-\theta_0/2 & {\mbox{ if $v(\xi)\leq-\theta_0/2$.}}
\end{matrix}
\right.$$
Then, $w$ is a
Lipschitz function in the
Heisenberg
group setting; so, by the Coarea Formula
in the
Heisenberg
group
(see~(1.4) in [MSC01]),
\begin{eqnarray*}
&&\int_{\BB_r(\xi_0)\cap \{ |v|\leq \theta_0/2\}}
|\gra v|=\\
&=&
\int_{\BB_r(\xi_0)}
|\gra w|=\\
&=& \int_{-\infty}^{+\infty} \,{\rm Per}_{\He^n} \Big(
\{ w=t\}\cap \BB_r(\xi_0)
\Big) \,dt\geq \\
&\geq&  \int_{-\theta_0/2}^{+\theta_0/2} \,{\rm Per}_{\He^n} \Big(
\{ v=t\}\cap \BB_r(\xi_0)
\Big) \,dt\,,\end{eqnarray*}
where ``${\rm Per}_{\He^n}$" denotes the
surface measure in~$\He^n$ (see~[MSC01] and~[FSSC01]
for details).

Consequently, using the
Isoperimetric Inequality in the
Heisenberg
group surface theory (see, e.g., Theorem~3.5
in~[FSSC01] and references
therein),
\begin{eqnarray*}
&& \int_{\BB_r(\xi_0)\cap \{ |v|\leq \theta_0/2\}}
|\gra v|\,\geq\,\\ &\geq&\,{\rm const}\,
\int_{-\theta_0/2}^{+\theta_0/2}
\Big(\min\Big\{
\Lgot^{2n+1}\Big(\BB_r(\xi_0)\cap \{v<t\}\Big)\,,\\ &&
\Lgot^{2n+1}\Big(\BB_r(\xi_0)\cap \{v>t\}\Big)
\Big\}\Big)^{(Q-1)/Q}\,dt\,\geq\\
&\geq&\,
{\rm const}\,
\Big(\min\Big\{
\Lgot^{2n+1}\Big(\BB_r(\xi_0)\cap \{v<-\theta_0/2\}\Big)\,,\\&&
\Lgot^{2n+1}\Big(\BB_r(\xi_0)\cap \{v>\theta_0/2\}\Big)
\Big\}\Big)^{(Q-1)/Q}
\,.\end{eqnarray*}
Therefore, using~\eqref{8uijnhjHHH}
and~\eqref{8uijnhjHHH2},
\begin{eqnarray*}
&&1+\int_{\BB_r(\xi_0)\cap \{ |u|\leq \theta_0\}}
|\gra u|\,\geq\,\\ &\geq&\,{\rm const}\,
\Big(\min\Big\{
\Lgot^{2n+1}\Big(\BB_r(\xi_0)\cap \{u<-\theta_0\}\Big)\,,\;\,\\&&
\Lgot^{2n+1}\Big(\BB_r(\xi_0)\cap \{u>\theta_0\}\Big)
\Big\}\Big)^{(Q-1)/Q}\,.\end{eqnarray*}
Accordingly,
\begin{equation}\label{9iJJJJJ}
\int_{\BB_r(\xi_0)\cap \{ |u|\leq \theta_0\}}
|\gra u|\,\geq\,
{\rm const}\,
r^{Q-1}\,,
\end{equation}
due to~\eqref{23ed5551}
and~\eqref{23ed5552}, provided that~$r_0$
is large enough.

We now fix~$A>0$ suitably big and
we use the Cauchy
Inequality and
Theorem~\ref{CC2}[{\rm i}],
to deduce from~\eqref{9iJJJJJ}
that
\begin{eqnarray*}
A \,\Lgot^{2n+1} \Big( \BB_r(\xi_0)\cap \{
|u|\leq\theta_0
\}\Big) +\frac{{{\,\rm const\,}}r^{Q-1}}{A}
\,\geq\\
\geq
\,A\, \Lgot^{2n+1} \Big( \BB_r(\xi_0)\cap \{
|u|\leq\theta_0
\}\Big) \,+\,
\frac1A\int_{\BB_r(\xi_0)}|\gra u|^2\,\geq\\
\geq \int_{\BB_r(\xi_0)\cap \{ |u|\leq \theta_0\}}
|\gra u|\,\geq\\
\geq\,{\rm const}\, r^{Q-1}\,,
\end{eqnarray*} from
which the claim in
Theorem~\ref{IMPROVE}
follows by taking~$A$ large enough.

\section{Proof of Theorem~\ref{classa-per}}\label{sec:5}

\subsection{Notions of periodicity in $\He^n$.} From
now on, we impose that~$F$ is intrinsically periodic, i.e., that it 
satisfies~\eqref{classperiodic}, and
we fix~$\omega=(\om_1,\dots,
\om_{2n})\in \Q^{2n}\setminus\{0\}$.
 
We now see how~$\omega$ naturally induces
a concept of periodicity for functions on the
Heisenberg group:

\begin{defn}\label{8jUU} $u:\He^n
\rightarrow \R$ is a function which is $\omega$-periodic if
$u(\xi)=u((k,0)\circ \xi)$ for any $k\in \Z^{2n}$ so that $\om\cdot k
=0$.
\end{defn}

We now suppose that
\begin{equation}\label{n3}
n\geq 2\,,\end{equation}
the case~$n=1$ being dealt with on page~\pageref{7ujhHHjJkK345}

We 
introduce the vectors $k^j\in\Z^{2n}$, 
for $j=1,\dots,2n$, in order to have
an integer base\label{6tgGGYTGG3ed} 
of the ${2n}$-dimensional
lattice generated by $\om$. For this, 
let us consider a rational orthogonal base of~$\Q^{2n}$ given
by the vectors~$\omega$,~$v^1$,~$\dots$, $v^{2n-1}$.

By~(\ref{n3}), it is easy to see\footnote{Indeed, we remark that ${\rm
Im}\,(\overline v^i v^j)$ is just the scalar product of $v^j$ with a
vector which we will denote 
$\bar v^i:=(-v^i_2,v^i_1,-v^i_4,v^i_3,\dots
-v^i_{2n},v^i_{2n-1})$.  
So, fixed~$k=1,\dots,2n-1$,
we use the base property to write $\bar v^k:=\sum_{i=1}^{2n-1} 
\alpha_i^k v^i+b^k\omega$, for suitable $\alpha_i^k,b^k\in\R$.
Thus, if 
$0={\rm
Im}\,(\overline v^i v^j)=\bar v^i\cdot v^j$ for all~$i$ and~$j$,
then $$\left\{\begin{array}{c} 0=\bar v^k\cdot
v^1=\alpha_1^k\,
|v^1|^2\\ \vdots\\ 0=\bar v^k\cdot
v^{2n-1}=\alpha_{2n-1}^k\,
|v^{2n-1}|^2\\ \end{array} \right. \,.$$ 
We 
have
thus
obtained that $\bar v^k=b^k\omega$ and, by the base property,
$b^k\neq0$.
But then one gets the contradiction
$$ 0\neq b^1 b^2 \,|\om|^2 = \bar v^1\cdot\bar v^2 =v^1 \cdot v^2=0 
\,.$$}
that there exists some~$i$ 
and~$j$ such that ${\rm
Im}\,(\overline v^i v^j)\,\neq\,0\,$. 
Let now~$Q\in \N$ be a
common multiple of the denominators of the coordinates of the vectors
$v^1,\dots, v^{2n-1},\omega$. Then, the vectors $k^j:= Qv^j$ for
$j=1,\dots, 2n-1$ and $k^{2n}:= Q\omega$ are the required integer base.
Furthermore, by
construction, 
\begin{equation}\label{8uhb7} \,{\rm Im}\,(\overline k^i
k^j)\,\neq\,0\,, \end{equation} for 
some~$i$ and~$j$ not equal to~$2n$.

Due to~\eqref{8uhb7}, and possibly exchanging~$k^i$ with~$k^j$,
we thus have that 
\begin{equation}\label{chow3}
\Theta_{ij}\,:=\,\,{\rm Im}\,(\overline k^i k^j)\,>\,0\,.
\end{equation}
If $m\in\R$, we define 
$k^{ij}_m:=m k^i+k^j$.
It follows from~(\ref{itekappa})
and~(\ref{chow3})
that
$$ (k^{ij}_m,0)\circ(-m k^i,0)\circ(-k^j,0)
\circ(z,t)=(z,t+2m\Theta_{ij})\,.
$$
Due to the above relation (taken with $m=1$) and~(\ref{chow3}),
an~$\omega$-periodic function is also periodic in 
the vertical variable, that is
\begin{equation} \label{chow4}
u(z,t)=u(z,t+2\Theta_{ij})\,.
\end{equation}
Of course, this holds true for any~$i$ and~$j$ such that $\Theta_{ij}
\neq 0$. 
Therefore, if~$u$ is~$\omega$-periodic, it is also~$\Theta$-periodic in 
the vertical direction, with~$\Theta$ being the greatest common divisors 
of the~$\Theta_{ij}\neq 0$.\label{7ujhn8ikj7ujhedfv}
Observe that, in particular,~$F(\cdot,u)$ 
is~$2$-periodic in the vertical 
direction.

\subsection{Constrained minimizers.}

We now fix~$M\geq 10$, \label{daqui}
to be chosen suitably large in the sequel. 
We seek the minimizers of~$\Fgot$ among the 
functions that are $\omega$-periodic and that 
satisfy a constraint on the 
hyperplanes $\{\xi\cdot (\omega,0)=-M |\om|\}$ and 
$\{\xi\cdot (\omega,0)=M|\om| \}$. 

We define
$$ U_\om\,:=\,\left\{ \sum_{j=1}^{2n-1}s_j k^j\,,\;\,s_j
\in \Big[-\frac 12,\frac 12\Big)
\right\}\,\times\,[-\Theta,\Theta)\,.$$
Geometrically, $U_\om$ is a $2n$-dimensional
``cube'' (or, better to say,
rectangular parallelopipedon) with edges of length $|k^1|,\dots,
|k^{2n-1}|, 2\Theta$. 

We now flow $U_\om$ in the $\om$-direction,
by using the 
Heisenberg group action. That is, we define
\begin{equation}\label{yhn98999}
S_{\om}\,:=\, \bigcup_{\alpha\in\R} (\alpha k^{2n},0)\circ 
U_\om
\,.\end{equation}
We also introduce the horizontal projection of $S_\omega$:
$$ H_\om \,:=\, \left\{ \sum_{j=1}^{2n-1}s_j k^j+\alpha k^{2n}
\,,\;\,s_j
\in \Big[-\frac 12,\frac 12\Big)\,,\,\;\alpha\in\R
\right\}\,.$$
Given $z=\sum_{1\leq j\leq 2n} z_j k^j\in \R^{2n}$, we define
$$ \zeta(z)\,:=\, 2\sum_{1\leq j\leq 2n}  z_j z_{2n}\,{\rm Im}\,
(\overline k^{2n} k^j)\,.$$
It is easily checked that~$S_\om$ may be seen as the union of intervals
of length~$2\Theta$ in the~$t$-direction, centered at
$\zeta(z)$, when $z$ varies in~$H_\om$, that is
\begin{equation}\label{7uHHHHyHGGG}
S_\om\,=\, \Big\{ 
(z,\,t+\zeta(z))\,,\,\; z\in H_\om\,,\,\; t\in[-\Theta,\Theta)
\Big\}\,.
\end{equation}

Fixed~$\delta\in (0,1/2)$, let\label{7yhnjj8}
\begin{equation}\label{consty}
\begin{array}{ll}
\cY_M\,:=\, &\Big\{u\in S^{1,2}_{\rm loc}\,, 
\;\,\omega-\mbox{periodic}, 
\\[1ex]
& u(\xi)\geq 1-\delta\ \mbox{for} \ \xi\cdot (\om,0)\geq 
M\,|\om|,\\[1ex]
 & u(\xi)\leq -1+\delta\ \mbox{for}\  
\xi\cdot (\om,0)\leq -M\,|\om|\Big\}\,.
 \end{array}\end{equation}

Then, we have the following existence result
for constrained minimizers:

\begin{prop} \label{ExistMinCon}
For any $M\geq 1$ 
there exists an absolute minimum in $\cY_M$ of
$$ \Fgot_{S_\om} (u)= \int_{S_\om} |\gra u(\xi)|^2+
F(\xi,u(\xi))\,d\xi\,.
$$
Moreover, if~$u^*$ is any absolute minimizer,
then
\begin{equation}\label{7ygvftgfFFF} 
\Fgot_{S_\om} (u^*)\,\leq \,c\,\Theta\, 
|k^1|\,\dots\,|k^{2n-1}|
\,,\end{equation}
where $c>0$ depends only on the structural constants
of $F$.
\end{prop}

\begin{proof} As usual, we use the notation~$\xi=(z,t)$
for points in~$\He^n$.
Set \begin{equation}\label{u0defuni}
u_0(\xi)\,:=\,
\left\{ 
\begin{matrix}
\frac{4\om}{|\om|}\cdot z & \,{\mbox{ if $|\frac{4\om}{|\om|}\cdot z
|\leq 1$,}}\\
1 & \,{\mbox{ if $\frac{4\om}{|\om|}\cdot z\geq 1$,}}\\
-1 & \,{\mbox{ if $\frac{4\om}{|\om|}\cdot z\leq -1$.}}\\
\end{matrix}
\right.
\end{equation}
It is easily seen that
$u_0\in \cY_M$.
 
Moreover, since~$k^{2n}$ is parallel to~$\om$, we conclude that
\begin{eqnarray*}
S_\om \cap \{|u_0|\neq 1\} &\subseteq& S_\om \cap\left\{
\Big|\frac{\om}{|\om|}\cdot z
\Big|\leq \frac 1 4
\right\}=\\ &=&
S_\om \cap\left\{
\Big|\frac{k^{2n}}{|k^{2n}|}\cdot z
\Big|\leq \frac 1 4
\right\}\subseteq\\ 
&\subseteq& 
\bigcup_{\alpha\in \left[-\frac{1}{4|k^{2n}|},
\frac{1}{4|k^{2n}|}\right]} 
(\alpha k^{2n},0)\circ 
U_\om
\,.
\end{eqnarray*}
Thence,
$$ \Fgot_{S_\om}(u_0)\,\leq\,c\,\Theta
|k^1|\,\dots\,|k^{2n-1}|
\,,$$
for a suitable~$c>0$.

Fix now~$R>0$ and let $$
D_R\,:=\, \Big\{\xi\, |\, \om\cdot z\in[-R|\om|,R|\om|]\Big\}
\cap S_\om\,.$$
Let also~$u_k$ be a minimizing sequence
for~$\Fgot_{S_\om}$.

By possibly cutting~$u_k$ at the $\pm 1$-levels,
which makes the energy decrease, we may assume
that \begin{equation}\label{ujm5678}
|u_k|\leq 1\,,\end{equation}
for any $k\in \N$.

Also,  
for~$k$ sufficiently large, $$ \Fgot_{S_\om}(u_k)\leq  
\Fgot_{S_\om}(u_0)<+\infty\,.$$ 
Hence,
$$ \sup_{k\in \N} \|\gra u_k\|_{L^2(S_\om)}\,<\,
+\infty\,.$$ 
This, (\ref{ujm5678}) and
standard embedding results (see, e.g.,~[VSCC92])
imply that there exists
a suitable function~$u$ so that,
up to
subsequences,~$\gra u_k$ weakly converges in~$L^2(D_R)$ 
to~$\gra u$ and~$u_k$ converges to~$u$ almost everywhere.
Therefore, the above mentioned weak convergence
gives that
$$ \liminf_{k\rightarrow +\infty} \int_{S_\om}|\gra u_k|^2
\geq
\liminf_{k\rightarrow +\infty} \int_{D_R}|\gra u_k|^2\geq
\int_{D_R}|\gra u|^2\,.$$
Fix now an arbitrarily small~$a>0$.
The pointwise convergence of~$u_k$,
Fatou's Lemma and the continuity of~$F$
when the second variable is in~$(-1,1)$
imply that 
$$ \liminf_{k\rightarrow +\infty} \int_{S_\om} F(\xi, u_k)\geq
\liminf_{k\rightarrow +\infty} \int_{{D_R}\atop{|u|\leq 1-a}} 
F(\xi, u_k)\geq
\int_{{D_R}\atop{|u|\leq 1-a}} F(\xi, u)\,.$$
Since~$a>0$ may be taken as small as we wish,
we conclude from the Monotone Convergence Theorem that
$$ \liminf_{k\rightarrow +\infty} \int_{S_\om} F(\xi, u_k)\geq
\int_{{D_R}} F(\xi, u)\,.$$
As a consequence,
$$ \inf_{\cY_M} \Fgot_{S_\om} 
=\liminf_{k\rightarrow+\infty}
\Fgot_{S_\om }(u_k)\geq \Fgot_{D_R} (u)\,.$$
Since~$R$ may be taken arbitrarily large, we
conclude the proof.
\end{proof} 

We now investigate further properties
of the constrained functions in the 
space~$\cY_M$.

\begin{lemma}\label{minmaxiny}
Let $u,\,v\in \cY_M$. Then,
$$ \min\{u,v\}\,,\;
\max\{u,v\}\in
\cY_M\,$$
and 
$$ \Fgot_{S_\om}(u)+\Fgot_{S_\om}(v)=
\Fgot_{S_\om}(\min \{ u,v\})+
\Fgot_{S_\om}(\max \{u,v\})\,,$$
provided that the above quantities are
finite.
\end{lemma} 

\begin{proof} 
%% We only prove the claim for the minimum,
%% because the case of the maximum 
%% follows by the same arguments.
%% Of course, if $\xi\cdot(\om,0)\leq -M|\om|$,
%% $$ \min\{u,v\}(\xi) \leq u(\xi)\leq -1+\delta\,.$$
%% On the other hand, if $\xi\cdot(\om,0)\geq M|\om|$,
%% then $u(\xi)\geq 1-\delta$ and $v(\xi)\geq 1-\delta$,
%% and so $$ \min\{u,v\}(\xi) \geq 1-\delta\,.$$
%% Finally, $$
%% \min\{u,v\}((k^i,0)\circ \xi) =\min
%% \{u((k^i,0)\circ \xi),\, v((k^i,0)\circ \xi)\}=$$
%% $$=\min\{u( \xi),\, v( \xi)\}=
%% \min\{u,\, v\}(\xi)\,,$$
%% for $i=1,\dots,2n-1$,
%% since~$u$ and~$v$ are $\om$-periodic, thence
%% $ \min\{u,v\}$ is $\om$-periodic too.
The first claim is obvious.
For the proof of the  second statement, just split
the domain of integration~${S_\om}$
into~$S_\om\cap \{u\le v\}$
and~$S_\om\cap\{ u> v\}$ and compute.
\end{proof}

We denote by~${\mathcal M}_M$
the set of constrained
minimizers given by Proposition~\ref{ExistMinCon}
(for a fixed~$\om$,
which is omitted in the notation).

\begin{lemma}\label{thJHJ}
If $u,\,v\in {\mathcal M}_M$,
then so do~$\min\{u,v\}$
and~$\max\{u,v\}$.
\end{lemma}

\begin{proof} 
Since
$u,\,v\in {\mathcal M}_M$,
$$ m:= \Fgot_{S_\om}(u)=\Fgot_{S_\om}(v)\leq
\Fgot_{S_\om}(w)\,,$$
for any $w\in \cY_M$.
But, by Lemma~\ref{minmaxiny}, both $ \min\{u,v\}$ and
$\max\{u,v\}$ are in
$\cY_M$ and therefore 
$$ \Fgot_{S_\om}(\min \{ u,v\})\geq m$$
and
$$ \Fgot_{S_\om}(\max \{u,v\})\geq m\,.$$
If either
$$ \Fgot_{S_\om}(\min \{ u,v\})> m$$
or
$$ \Fgot_{S_\om}(\max \{u,v\})>m\,,$$
we would have that
$$ \Fgot_{S_\om}(u)+
\Fgot_{S_\om}(v)=
2m <
\Fgot_{S_\om}(\min \{ u,v\})+
\Fgot_{S_\om}(\max \{u,v\})\,,$$
which is in contradiction with 
Lemma~\ref{minmaxiny}. Therefore,
$$
\Fgot_{S_\om}(\min \{ u,v\})=
\Fgot_{S_\om}(\max \{u,v\})=m\,,$$
proving that
$\min\{u,v\}$
and~$\max\{u,v\}$ are in~${\mathcal M}_M$.
\end{proof}

\begin{lemma}\label{ascolile}
Let $v_n\in {\mathcal M}_M$, for any $n\in\N$.
Then, there exists a subsequence $n_k\rightarrow +\infty$
and a function $v\in{\mathcal M}_M$
such that~$v_{n_k}$
converges to~$v$ uniformly on compact sets.
\end{lemma}

\begin{proof} By [M95],~$v_n$ is an equicontinuous
family. Also, $|v_n|\leq 1$, 
and thus~$v_n$ converges
locally uniformly
to some~$v$, up to subsequences, by the Theorem of Ascoli.
Then, by Fatou's Lemma,
\begin{equation}\label{thnYY}
\liminf_{n\rightarrow+\infty} 
\int_{S_\om} F(\xi,v_n(\xi))\,d\xi\,\geq\,
\int_{S_\om} F(\xi,v(\xi))\,d\xi\,.
\end{equation}

Since $\cY_M$ is obviously closed under
pointwise limit and~$v_n\in\cY_M$,
we also have that
\begin{equation}\label{thnYY2}
v\in \cY_M\,.
\end{equation}

Furthermore, since all the~$v_n$'s are minimizers,
it follows that $\Fgot_{S_\om} (v_n)$ takes 
always the
same value for any $n\in\N$, say
\begin{equation} \label{ujKK}
\Fgot_{S_\om} (v_n)=:m\geq 0\,.\end{equation}
In particular,
$$ \int_{S_\om} |\gra v_n|^2 \leq m\,,$$
thus, up to subsequences, we may
assume that~$\gra v_n$ converges to $\gra v$
weakly in~$L^2$.
Therefore,
$$
\liminf_{n\rightarrow+\infty} 
\int_{S_\om} |\gra v_n(\xi)|^2\,d\xi\,\geq\,
\int_{S_\om} |\gra v(\xi)|^2\,d\xi\,.
$$ Combining this with~(\ref{thnYY})
and~(\ref{ujKK}),
we thus get that $\Fgot_{S_\om} (v)\leq m$. Therefore,
by~(\ref{thnYY2}), we conclude 
that~$v\in{\mathcal M}_M$. 
\end{proof}

We now define the minimal
minimizer, that is, the pointwise
infimum of all the minimizers in~${\mathcal M}_M$:
\begin{equation} \label{defmimi}
u_M (\xi)\,:=\,\inf_{u\in {{\mathcal M}_M}}
u(\xi)\,,\end{equation}
for any $\xi\in\He^n$. 
Note that $u_M$ is indeed a minimizer, according
to the following observation:

\begin{lemma}\label{yhnUUjJ87}
$u_M\in{{\mathcal M}_M}$.
\end{lemma}

\begin{proof} For any $\xi\in\He^n$, we consider
a sequence $u_{n,\xi}\in {\mathcal M}_M$
in such a way
\begin{equation}\label{yHHHjJJkH}
\lim_{n\rightarrow +\infty} u_{n,\xi} (\xi)=
\inf_{u\in{\mathcal M}_M} u(\xi)=u_M(\xi)\,.
\end{equation}
Let us now write the countable set $\Q^{2n+1}$
as $\{ \eta_1,\eta_2, \dots\}$.
We define
\begin{equation} \label{YHJ7yy667}
v_n(\xi)\,:=\, \min\{
u_{n,\eta_1}(\xi)\, , \, \dots \, , \,
u_{n,\eta_n}(\xi) 
\}\,.
\end{equation} From Lemma~\ref{thJHJ},
we conclude that~$v_n\in {\mathcal{M}}_M$
for any~$n\in\N$. Thus, by
Lemma~\ref{ascolile},
perhaps passing to a subsequence, we have that
\begin{equation}\label{ghghghghhh}
{\mbox{
$v_n$ 
pointwise converges to a suitable $v$}}
\end{equation}
and that
\begin{equation}\label{uYYhG88}
v\in {\mathcal M}_M\,.
\end{equation}
We show that
\begin{equation}\label{ynJJJ}
u_M \,\geq\, v\,.
\end{equation}
The proof of~(\ref{ynJJJ}) 
is by contradiction.
If~(\ref{ynJJJ}) were false, there would
exist~$\e_0>0$ and $\tilde\xi\in\He^n$
such that
$$ u_M (\tilde\xi)+6 \e_0
\leq v(\tilde\xi)\,.$$
In the light of~(\ref{defmimi}),
there exist~$\tilde u\in {\mathcal M}_M$
such that
$$ \tilde u(\tilde\xi) \leq u_M (\tilde\xi)+\e_0$$
and therefore
$$ \tilde u
(\tilde\xi)+5 \e_0
\leq v(\tilde\xi)\,.$$
Since both~$v$ and~$\tilde u$ are in ${\mathcal M}_M$,
their modulus of continuity can be
bounded uniformly, thanks to~[M95].
Thus, since~$\Q^{2n+1}$ is dense in~$\He^n$,
there exists~$\tilde m\in \N$ in such a way
\begin{eqnarray*}
&&|v(\tilde\xi)-v(\eta_{\tilde m})|\leq \e_0
\,{\mbox{ and}}\\
&&|\tilde u(\tilde\xi)-
\tilde u(\eta_{\tilde m})|\leq \e_0
\,.
\end{eqnarray*}
Therefore,
\begin{equation}\label{tghJJJhFFR} 
\tilde u
(\eta_{\tilde m})+3 \e_0
\leq v(\eta_{\tilde m})\,.
\end{equation} 
Recalling~(\ref{yHHHjJJkH})
and~(\ref{ghghghghhh}), we now
take $\tilde n\in\N$ in such a way
\begin{equation}\label{yhJKJK76544}
\tilde n\geq \tilde m\,,\end{equation}
\begin{equation} \label{ynJII7666}
u_{\tilde n, \eta_{\tilde m}}(\eta_{\tilde m})
\leq \inf_{u\in {\mathcal M}_M} 
u(\eta_{\tilde m})+\e_0
\end{equation}
and
\begin{equation}\label{766TTyT}
v(\eta_{\tilde m}) \leq 
v_{\tilde n} (\eta_{\tilde m})+\e_0\,.
\end{equation}
In particular, from~(\ref{ynJII7666}),
$$ u_{\tilde n, \eta_{\tilde m}}(\eta_{\tilde m})
\leq
\tilde u(\eta_{\tilde m})
+\e_0$$
and therefore, by~(\ref{tghJJJhFFR}),
(\ref{766TTyT}), (\ref{yhJKJK76544})
and~(\ref{YHJ7yy667}),
$$
u_{\tilde n, \eta_{\tilde m}}(\eta_{\tilde m})
+\e_0 \leq
v_{\tilde n} (\eta_{\tilde m})\leq 
u_{\tilde n, \eta_{\tilde m}}
(\eta_{\tilde m})
\,,
$$
which is a contradiction 
and thus gives the proof of~(\ref{ynJJJ}).

In fact, by~(\ref{defmimi})
and~(\ref{uYYhG88}), we conclude that $u_M\leq v$.
Therefore, by~(\ref{uYYhG88}) and~(\ref{ynJJJ}),
$$ u_M=v\in {\mathcal M}_M\,,$$
as desired.
\end{proof}

We also point out that, since~$u_M$ is a
minimizer (due to Lemma~\ref{yhnUUjJ87}),
by cutting~$u_M$ at the levels~$\pm 1$,
it easily follows that~$|u_M|\leq 1$.

We now observe that the minimal minimizer
inherits the vertical periodicity
of the functional, thus improving the one in~(\ref{chow4}).
This fact will be of crucial importance in the proof of
Lemma~\ref{centrle} later on.

\begin{lemma}\label{2periodicity}
$u_M$ is $2$-periodic in the vertical variable.
\end{lemma}

\begin{proof} If $\Theta=1$ there is nothing to prove, otherwise let
$u_{M,j}(z,t):=u_M(z,t+2j)$. Since~$F$ is~$2$-periodic in the vertical
variable
and $\gra $ is invariant with respect to vertical
translations, then~$u_{M,j}$ is also a minimizer for any
$j=0,\dots,\Theta-1$. Hence $v(\xi)=\min\{u_M(\xi),u_{M,1}(\xi),\dots
,u_{M,\Theta-1}(\xi)\}$ is also a 
minimizer by Lemma~\ref{thJHJ}. Hence,
by the definitions of $u_M$ and $v$, 
$$u_M\leq v\leq u_M$$ and so~$u_M=v$.
But~$v$ is~$2$-periodic in~$t$
and thus so is~$u_M$. \end{proof}

\begin{rem} Due to Lemma~\ref{2periodicity}, without 
loss of generality, one may suppose that 
the functions of~$\cY_M$ 
are~$2$-periodic in the vertical variable, that $\Theta=1$ 
and that 
$$ U_\om\,:=\,\left\{ \sum_{j=1}^{2n-1}s_j k^j\,,\;\,s_j
\in \Big[-\frac 12,\frac 12\Big)
\right\}\,\times\,[-1,1)\,.$$
\end{rem}

\subsection{A group action on the space of functions.}
Given a function $u:\He^n\longrightarrow
\R$ and $k\in \Z^{2n}$, we define
$$ T_k u(\xi)\,:=\, u\Big( (k,0)\circ
\xi\Big)\,,
$$
for any $\xi\in \He^n$.

\begin{lemma}\label{yhnUU5tefg}
Let $k\in\Z^{2n}$. If $u$ is $\om$-periodic and $2$-periodic in the 
vertical direction,
then so is $T_k u$.
\end{lemma}

\begin{proof} Let $h\in \Z^{2n}$, with $\om\cdot h=0$. 
Let $a,b\in\Z$ be
so that
\begin{eqnarray*}
(h,0)\circ(k,0)&=&(h+k,2a)\qquad\;{\mbox{and}}\\ 
(k,0)\circ(h,0)&=&(h+k,2b)\,.\end{eqnarray*}
By the periodicity of~$u$,
$$T_k u \Big((h,0)\circ \xi\Big)=u(h+k+z,t+2b)=u(h+k+z,t+2a)=$$
$$=u\Big( (h,0)\circ (k,0)\circ\xi\Big)=u\Big( (k,0)\circ\xi
\Big)=
T_k u(\xi)\,.$$
Thus,~$u$ is $\om$-periodic. 

The vertical $2$-periodicity is obvious.
\end{proof}

We now show that~$\cY_M$ behaves
nicely with respect to the action~$T_k$:

\begin{lemma}\label{minmaxintras}
Let $k\in \Z^{2n}$ be so that
$\om \cdot k\geq 0$. Let
$u\in\cY_M$ be $2$-periodic in the vertical
direction. Then,
$$ \min\{u,\,T_k u\}\,\in\,
\cY_M\,.$$
Analogously, if $k\in \Z^{2n}$ is
so that
$\om \cdot k\leq 0$ and
$u\in\cY_M$, then
$$ \max \{u,\,T_k u\}\,\in\,
\cY_M\,.$$
\end{lemma}

\begin{proof} 
We prove the first claim, the second one being
analogous.
If $\xi\cdot(\om,0)\leq -M|\om|$,
$$ \min\{u,\,T_k u\}(\xi) \leq u(\xi)\leq -1+
\delta\,,$$
because $u\in\cY_M$.

If
$\xi\cdot(\om,0)\geq M|\om|$,
let $\xi':= (k,0)\circ \xi$.
Then,
$$ \xi' \cdot (\om ,0)=\xi\cdot (\om,0)+
\om\cdot k\geq M\,|\om|+0=M\,|\om|\,,$$
and so
$u(\xi)\geq 1-\delta$
and $u(\xi')\geq 1-\delta$, since
$u\in\cY_M$.
Therefore,
$$ \min\{u,\,T_k u\}(\xi)=
\min\{u(\xi)\,,u(\xi') \} \geq 1-\delta\,.$$
Finally, by Lemma~\ref{yhnUU5tefg},
$ \min\{u,\,T_k u\}$
is $\om$-periodic.
\end{proof}

%%% The group action~$T_k$
%%% plays a crucial role when our functional
%%% is~$\Z^{2n}$-periodic with respect to the
%%% Heisenberg group action.
Given $S\subset \He^n$, we define
\begin{eqnarray*}
T_k S &:=& \{ (-k,0)\circ \xi,\;\,\xi
\in S  \}\quad\quad{\mbox{and}}\\
T_k\cY_M &:=& \{ T_k u,\;\,u
\in \cY_M\}\,.
\end{eqnarray*}
Analogously, we set
$$ \Fgot_{T_k S_\om}(u):=
\int_{T_k S_\om} |\gra u(\xi)|^2+F(\xi,u(\xi))\,d\xi\,.$$
Then, since the functional~$\Fgot$ is $\Z^{2n}$-periodic 
with respect to the
Heisenberg group action, then it
is
invariant under the action of~$T_k$, according
to the following observation, whose elementary
proof is omitted:
 
\begin{lemma}\label{invfunztra}
We have that
$$ \Fgot_{T_k S_\om} (T_k u) = \Fgot_{S_\om} (u)\,,$$
for any $k\in\Z^{2n}$.
\end{lemma}

%% \begin{proof} It is an immediate consequence of 
%% the fact that the vector fields are invariant 
%% with respect to the left groupa ction and 
%% that the Lebesgue measure is the Haar measure of 
%% the group.
%% \end{proof}

We will see in Lemma~\ref{yuuiiu}
here below that
Lemma~\ref{invfunztra} may in fact
be straightened for $\om$-periodic functions.
For this scope, we need to better investigate
the invariance properties of~$S_\om$. First of all,
it is obvious from~\eqref{yhn98999}
that
\begin{equation}
T_{\beta k^{2n}} S_\om \,=\,S_\om\,,
\end{equation}
for any $\beta\in\R$.
Moreover,

\begin{lemma}\label{yuuiiu}
Let~$u$ be $\om$-periodic.
Then,
$\Fgot_{S_\om}(u)=\Fgot_{T_v S_\om} (u)$, for any~$v\in
\R^{2n}$.
\end{lemma}
\begin{proof} 
The result is a consequence of the following 

{\em Claim:} Let $v\in\R^{2n}$ and $f$ be an $\om$-periodic, integrable function.
Then, 
$$ \int_{T_v S_\om} f\,=\,
 \int_{S_\om} f\,.$$

\bigskip
We observe that, for any fixed $z\in\R^{2n}$,
the map $t\mapsto f(z,t)$ is $2\Theta$-periodic,
according to~\eqref{chow4}, and therefore
\begin{equation}\label{8ijmnJNHHJH}
\int_{-\Theta}^{\Theta} 
f(z,t+\tau)\,dt\,=\,\int_{-\Theta}^{\Theta}f(z,t)\,dt\,,
\end{equation}
for any $\tau\in\R$ and any
$z\in\R^{2n}$.

Let also 
$$ g(z)\,:=\,\int_{-\Theta}^{\Theta}f(z,t)\,dt\,.$$
Exploiting~\eqref{8ijmnJNHHJH}
and the fact that~$f$ is~$\om$-periodic, it follows that
\begin{eqnarray*}
g(z+k^j)&=& \int_{-\Theta}^{\Theta}f(z+k^j,t)\,dt=\\
&=& \int_{-\Theta}^{\Theta}f\Big((k^j,0)\circ (z,t)\Big)\,dt
= \\&=&\int_{-\Theta}^{\Theta}f(z,t)\,dt =\\ &=& g(z)
\,,\end{eqnarray*}
for $j=1,\dots, 2n-1$.

Consequently, 
$$ \int_{H_\om} g(z-v)\,dz
\,=\,
\int_{H_\om} g(z)\,dz\,,$$
for any $v\in\R^{2n}$.

Therefore, recalling also~\eqref{7uHHHHyHGGG},
we conclude that
\begin{eqnarray*}
\int_{T_v S_\om} f &=& \int_{S_\om} f\Big((-v,0)\circ \xi\Big)\,d\xi
\,=
\\ &=& \int_{H_\om}\int_{-\Theta+\zeta(z)}^{\Theta+\zeta(z)} 
f\Big(z-v, 
t+2\,{\rm 
Im}\,(\overline v z) \Big) \,dt\,dz\,=\\
&=& \int_{H_\om}\int_{-\Theta}^\Theta f (z-v, t
) \,dt\,dz\,=\\ &=&
\int_{H_\om}
g(z-v)
\,dz\,=\\ &=& 
\int_{H_\om}
g(z)
\,dz\,=\\ &=&
\int_{S_\om} f \,,
\end{eqnarray*}
as desired.
\end{proof}

In analogy with Proposition~\ref{ExistMinCon},
we denote by $T_k{\mathcal M}_M$ the
set of minimizers of the functional~$\Fgot_{T_k S_\om}$
on the space~$T_k \cY_M$.

\begin{lemma}\label{7uhjn9ik}
We have that
$u\in {\mathcal M}_M$ if and only if
$T_k u\in T_k {\mathcal M}_M$, for any $k\in\Z^{2n}$.
\end{lemma}

\begin{proof} We prove the ``only if'' part.
Let~$w\in T_k\cY_M$.
Then, $w=T_k v$, for some~$v\in \cY_M$.
Thus, exploting that $u\in  
{\mathcal M}_M$ and Lemma~\ref{invfunztra},
we conclude that
$$ \Fgot_{T_k S_\om}(T_k u)=\Fgot_{S_\om}(u)\leq
\Fgot_{S_\om}(v)=\Fgot_{T_k S_\om}({T_k v})=
\Fgot_{T_k S_\om}(w)\,,$$
thence $T_k u$ is a minimizer of~$\Fgot_{T_k S_\om}$
in~$T_k \cY_M$.
\end{proof}

\begin{lemma}\label{dopo0.10}
Let $k\in\Z^{2n}$
such that $\om\cdot k\geq0$. Let $u\in 
{\mathcal M}_M$ be $2$-periodic in the
vertical direction. Then,
$\min \{u,\,T_k u\}\in
{\mathcal M}_M$.
\end{lemma}

\begin{proof}
This is a variation of the proof of Lemma~\ref{thJHJ}.
We provide full details for the facility of the reader.
Let $m:=\Fgot_{S_\om}(u)$. We show that
\begin{equation}\label{qqa}
\Fgot_{S_\om}\Big( 
\min \{u,\,T_k u\}
\Big)\,=\,m\,.
\end{equation}
To prove~(\ref{qqa}) suppose, by contradiction,
that
\begin{equation}\label{qqb}
\Fgot_{S_\om}\Big(
\min \{u,\,T_k u\}
\Big)\,\not=\,m\,.
\end{equation}
By Lemma~\ref{minmaxintras},
we know that 
\begin{equation}\label{qqx}
\min \{u,\,T_k u\}\in \cY_M\,.\end{equation}
Thus,~(\ref{qqb})
gives that
\begin{equation}\label{qqc}
\Fgot_{S_\om}\Big(
\min \{u,\,T_k u\}
\Big)\, > \,m\,.
\end{equation}
Let $h:=-k$. By Lemma~\ref{minmaxintras},
$$ \max \{u,\,T_h u\}
\in \cY_M$$
and so
$$ \max  \{T_k u,\, u\}=
T_k\Big(
\max \{u,\,T_h u\}
\Big)\in T_k \cY_M\,.
$$ Consequently, from
Lemma~\ref{7uhjn9ik},
$$ \Fgot_{T_k S_\om}(T_k u)
\,\leq\,
\Fgot_{T_k S_\om} \Big( 
\max
\{u,\,T_k u\}
\Big)\,.$$ This and Lemma~\ref{invfunztra}
give that
$$ m\,\leq\,
\Fgot_{T_k S_\om} \Big(
\max
\{u,\,T_k u\}
\Big)\,.$$ 
Therefore, by Lemma~\ref{yuuiiu},
\begin{equation}\label{qqd}
m\,\leq\,
\Fgot_{S_\om} \Big(
\max
\{u,\,T_k u\}
\Big)\,.
\end{equation}
Moreover, using again Lemmata~\ref{invfunztra} 
and~\ref{yuuiiu},
\begin{equation}
\Fgot_{S_\om} (T_k u)=\Fgot_{T_k S_\om}(u)=
\Fgot_{S_\om}(u)=m\,.
\end{equation}
Consequently, by~(\ref{qqc}),~(\ref{qqd})
and Lemma~\ref{minmaxiny},
$$2m =\Fgot_{S_\om}(u)+
\Fgot_{S_\om}(T_k u)=
$$ 
$$=\Fgot_{S_\om}\Big( \min
\{u,\,T_k u\}
\Big)+
\Fgot_{S_\om}
\Big(\max\{u,\,T_k u\}
\Big)>2m\,.$$
This contradiction gives the proof
of~(\ref{qqa}).

The desired result
thus follows from~(\ref{qqa}) and~(\ref{qqx}).
\end{proof}

\subsection{The Birkhoff property.}
In Mather theory, it is quite common
to seek orbits which possess some
kind of monotonicity, usually referred to
as Birkhoff property. Such property has then
been extended and used also in~[CL98], [CL01] and~[RS03]. 

We now adapt the concept of
Birkhoff property, in order to deal with the 
Heisenberg group. 

\begin{defn}\label{UJNgJhjhg} 
A function $u:\He^n\longrightarrow\R$
is said to satisfy the Birkhoff
property with respect to $\om$ if
for any $k\in \Z^{2n}$ so that
$\om\cdot k\geq 0$, we have that
$T_k u(\xi)\geq u(\xi)$ for any $\xi\in \He^n$.
\end{defn}

We observe that, if~$u$ satisfies the Birkhoff
property and~$k\in\Z^{2n}$ is such that
$\om\cdot k\leq 0$, then, writing $k':=-k$ and
$\xi':= (k,0)\circ \xi$, we have
$\om\cdot k'\geq 0$ and thus
\begin{equation}\label{ujnKKl} u(\xi) =
T_{k'} (\xi')\geq u(\xi')=T_k u(\xi)\,,
\end{equation}
for any $\xi\in \He^n$.

\begin{lemma}\label{4rgbu}
Suppose that~$u$ verifies the Birkhoff
property with respect to~$\om$.
Let $\theta\in\R$ and assume that
\begin{equation}\label{ujnKKl2}
S\subseteq \{ u< \theta \}\,.
\end{equation}
Then,
$$ T_{k} S\subseteq \{ u< \theta \}\,,$$
for any $k\in\Z^{2n}$ so that
$\om\cdot k\geq0$. 
\end{lemma}

\begin{proof} Take $k$ as above and
$\xi\in T_k S$. Then, there exists $\eta\in S$
so that $\xi=(h,0)\circ \eta$, with $h:=-k$.
Note that $\om\cdot h \leq 0$, thus,
exploiting~(\ref{ujnKKl}) and~(\ref{ujnKKl2}),
we conclude that
$$ \theta > u(\eta)\geq T_h u(\eta)=
u((h,0)\circ \eta)=u(\xi)\,,$$
as desired.
\end{proof}

We show that the minimal minimizer $u_M$
is Birkhoff (with respect to the vector
$\om$ in the definition of $S_\om$):

\begin{lemma}\label{uMbir}
$u_M$ satisfies the Birkhoff property.
\end{lemma}

\begin{proof}
Let $k\in \Z^{2n}$ such that
$\om \cdot k\geq 0$. Then,
by Lemmata~\ref{yhnUUjJ87}, \ref{2periodicity}
and~\ref{dopo0.10},
$$ \min \{ u_M,
\, T_k u_M\}\,\in\,{\mathcal M}_M\,.$$
Therefore, by~(\ref{defmimi}),
$$u_M \,\leq\,
\min \{ u_M,
\, T_k u_M\}\,\leq\,T_k u_M\,,$$
as desired.\end{proof}

Here is an interesting geometric
property of Birkhoff functions:

\begin{lemma}\label{centrle}
Let~$\theta \in\R$, $a,r>0$.
Suppose that~$u$ is~$2$-periodic in the vertical
direction, that it
satisfies the Birkhoff property 
with respect to~$\om$ and that
\begin{equation}\label{7gfdrR}
\BB_r(\xi_0) \subseteq \{ u< \theta \}\,.\end{equation}
Then, for any $\xi\in \He^n$ such that
$(\om,0)\cdot (\xi-\xi_0)
\leq -a|\om|$, we have that $u(\xi)<\theta$,
as long as~$r$ and~$a$ 
are larger than a suitable positive universal constant.
\end{lemma}

This result somehow 
guaranties
that if, 
say~$u>1-\delta$ in a ball, then it is so in a half space. 

\begin{proof}
Let $\xi=(z,t)$ as requested here above.
Denote also $\xi_0=(z_0,t_0)$.
Let us consider the tiling~$\Tgot^*$
of $\R^{2n}$ made by standard Euclidean
cubes of side~$1$.
More explicitly,~$\Tgot^*$
is the collection of the cubes
$K+[-1/2,1/2)^{2n}$, where~$K$ varies in~$\Z^{2n}$.
Let $Q^*$ be the cube
of~$\Tgot^*$ that contains $z_0$. 
Then, there exists
$p\in\Z^{2n}$ so that
\begin{equation}\label{7uhnjj}
z\in p+Q^*\,.\end{equation} Then,
\begin{equation}\label{8iujnokm}
\frac{\omega }{|\om|}\cdot p \leq -a+\sqrt{2n}\leq 0\,,
\end{equation}
provided that $a$ is large enough.

Also, by \eqref{7uhnjj}, we get that there exists $w\in Q^*$
such that $z=p+w$. Then,
there exists $\tau\in\R$ so that 
\begin{equation}\label{345678o}
(p,0)\circ (w,T)=(z,T+\tau)\,,\end{equation}
for any $T\in\R$. 
Thus, let $m$ be the unique integer so that
$$ \frac{t-t_0-\tau-1}{2}\leq m<
\frac{t-t_0-\tau+1}{2}$$
and define $T:=t-\tau-2m$.
Then, $|T-t_0|\leq 1$ and so, if~$r$
is large enough, $(w,T)\in \BB_r(\xi_0)$.

Moreover, by~\eqref{8iujnokm}, \eqref{7gfdrR} and
Lemma~\ref{4rgbu},
$$ (p,0)\circ (w,T)\,\in\, 
T_{-p}\Big(\BB_r (\xi_0)\Big)
\,\subseteq\, \{ u<\theta\}
\,.$$
Accordingly, by the vertical periodicity
and~\eqref{345678o},
\begin{eqnarray*}
&& u(\xi)=u(z,T+\tau+2m)=u(z,T+\tau)=\\
&& = u\Big( (p,0)\circ (w,T)\Big)<\theta\,,
\end{eqnarray*}
as desired.
\end{proof} 

\subsection{Koranyi balls outside the interface.}

Let us recall that in Theorem~\ref{CC2} we introduce 
a radius~$r_0$ which depends only on the 
structural data.  
We now show that we can get a Koranyi ball of universally 
large radius, say~$r_0$,
not touching the interface.

For this, we first need to enlarge our domain $S_\om$
so that the balls of radius~$r_0$ may comfortably fit in. 
We also need to enlarge our periodicity
mesh, in order to
obtain local minimizers on 
page~\pageref{8iujJJjUuHhUY5678}.

To these effects, given $p\in\N$, we define
\begin{eqnarray}\label{eqnuuHHy}
{\mathcal{Y}}_M^p &:=& \Big\{
u\in S^{1,2}_{\rm loc}\,,\quad u\Big( (pk,0\Big)\circ \xi
\Big)=u(\xi)=u\Big( (0,2p\Theta)\circ\xi\Big)\,,\nonumber\\
&& {\mbox{for any }} \xi\in\He^n{\mbox{ and }}
k\in\Z^{2n}{\mbox{ s.t. }}\om\cdot k=0\,,\\
&& u(\xi)\geq 1-\delta {\mbox{ for }} \xi\cdot (\om,0)
\geq M\,|\om|\,,\nonumber\\
&& u(\xi)\leq -1+\delta {\mbox{ for }} \xi\cdot (\om,0)
\leq -M\,|\om|\,\Big\}\nonumber
\,.\end{eqnarray}
Note that ${{\cY}}_M^1={{\cY}}_M$, 
due to~(\ref{chow4})
and~(\ref{consty}).

Given $v\in \R^{2n+1}$ and $S\subset \He^n$, we define
$$ T_v S \,:=\, \{ (-v)\circ \xi\,,\quad\xi\in S\}\,.$$

We also use the following notation. Given
$\ell\in\Z^{2n}$
with $\om\cdot \ell=0$, we write
$$ \ell = \sum_{1\leq j\leq 2n-1}\ell^j k^j\,.$$
Then, for $p\geq 1$, we define
\begin{equation}\label{equasdfgh5tg}
S_\om^p \,:=\,
\bigcup_{
{{\ell\in\Z^{2n},\;
\om\cdot \ell=0
}\atop{0\leq \ell^1,\dots\ell^{2n-1}
\leq p-1}}\atop{0\leq m\leq p-1}
} T_{(0,2 m\Theta)}T_{(\ell,0)} S_\om
\,.\end{equation}
Of course, $S_\om^1=S_\om$.
Also, if~$m\in\Z$ and~$\ell\in\Z^{2n}$ in such a way
$\om\cdot \ell=0$,
it can be seen from~(\ref{7uHHHHyHGGG})
that
\begin{equation}\label{6yhYh8ujm}
T_{(0,2 m\Theta)}T_{(\ell,0)} S_\om\,=\,
\Big\{
(z,t-2m\Theta+\zeta_\ell (z))\,,\quad z\in H_{\om \ell}\,,
t\in [-\Theta,\Theta)
\Big\}\,,
\end{equation} 
with
\begin{eqnarray*}
H_{\om \ell} &:=& \Big\{ 
z\in \He^n \;\,|\,\; z+\ell\in H_\om
\Big\}\quad\;{\mbox{and}}\\
\zeta_\ell (z)&:=& \zeta(z+\ell) -2\,{\rm Im}\, (\overline \ell z)\,.
\end{eqnarray*}
In particular, it follows that
\begin{equation}\label{dolla}{\mbox{
the union 
in~(\ref{equasdfgh5tg})
is non-overlapping. 
}}\end{equation}

We now observe that~$S_\om^p$ enjoys a good
translation invariance property:

\begin{lemma}\label{7ujhHhyhH}
Let~$u\in \cY_M^p$. Then,
$\Fgot_{S_\om^p}(u)=\Fgot_{T_v S_\om^p}(u)$
for any $v\in\R^{2n+1}$.\end{lemma}

\begin{proof}
The proof is a modification of
the one of Lemma~\ref{yuuiiu}.
The fact that the Heisenberg group elements do not
commute makes the argument quite technical,
thence we provide full details for the reader's 
convenience.

The claim in Lemma~\ref{7ujhHhyhH}
follows from the following one:
if~$f$ is an integrable function
so that
\begin{equation}\label{6yhHhHhJj7ujI90} 
f\Big( (pk,0)\circ \xi
\Big)=f(\xi)=f\Big( (0,2p\Theta)\circ\xi\Big)
\end{equation}
for any $\xi\in\He^n$ and any
$k\in\Z^{2n}$ so that $\om\cdot k=0$, 
then
\begin{equation}\label{7uhHhHujJjj6789}
\int_{T_v S_\om^p} f\,=\, \int_{S_\om^p} f\,,
\end{equation}
for any~$v\in \R^{2n+1}$.

In order to prove~(\ref{7uhHhHujJjj6789}),
we deduce from~(\ref{6yhHhHhJj7ujI90})
that
$$ f(z,t+2\Theta\,p)=f\Big( (0,2\Theta p)\circ (z,t)\Big)
=f(z,t)$$ and so,
for any fixed $z\in\R^{2n}$,
the map $t\mapsto f(z,t)$ is $2\Theta p$-periodic.

Accordingly,
\begin{equation}\label{8ijmnJNHHJH-bis}
\int_{-p\Theta}^{p\Theta}
f(z,t+\tau)\,dt\,=\,\int_{-p\Theta}^{p\Theta}f(z,t)\,dt\,,
\end{equation}
for any $\tau\in\R$ and any
$z\in\R^{2n}$.

We now define
$$ g(z)\,:=\,\int_{-p\Theta}^{p\Theta}f(z,t)\,dt\,.$$
It follows from~(\ref{6yhHhHhJj7ujI90})
and~(\ref{8ijmnJNHHJH-bis})
that
\begin{eqnarray*}
g(z+pk^j)&=& \int_{-p\Theta}^{p\Theta}f\Big(
(pk^j, 0)\circ (z,t-2p\,{\rm Im}\, (\overline k^j z)
\Big)\,dt=\\
&=& \int_{-p\Theta}^{p\Theta}f\Big(
z,t-2p\,{\rm Im}\, (\overline k^j z)
\Big)\,dt
= \\&=&\int_{-p\Theta}^{p\Theta}f(z,t)\,dt =\\ &=& g(z)
\,,\end{eqnarray*}
for any $1\leq j\leq 2n-1$.

Then, if
\begin{eqnarray*} H_\om^p&:=& \bigcup_{
{\ell\in\Z^{2n},\;\om\cdot \ell=0}\atop{0\leq \ell^{j}\leq p-1}
}
H_{\om \ell}\,=\\ &=& \left\{ 
\sum_{j=1}^{2n-1} s_j k^j+\alpha k^{2n}\,,\;\, s_j\in \Big[
-\frac 1 2,\frac 12 +p-1
\Big]\,,\;\alpha\in\R
\right\}\,,\end{eqnarray*}
we conclude that
\begin{equation}\label{ujJjJJj}
\int_{H^p_\om}g(z-w)\,dz\,=\,
\int_{H_\om^p} g(z)\,dz\,,
\end{equation}
for any $w\in\R^{2n}$.

Take now any $v=(v',v_{2n+1})\in\R^{2n+1}$.
We exploit~(\ref{equasdfgh5tg}),
(\ref{6yhYh8ujm}), (\ref{dolla})
(\ref{8ijmnJNHHJH-bis}) and~(\ref{ujJjJJj})
to obtain that
\begin{eqnarray*}
&& \int_{S_\om ^p} f(v\circ \xi)\,d\xi =\\
&=& \sum_{
{{\ell\in\Z^{2n},\;\om\cdot \ell=0
}\atop{0\leq \ell^{j}
\leq p-1}}\atop{0\leq m\leq p-1}
} \int_{H_{\om \ell}}\int_{\zeta_\ell(z) -(2m+1)\Theta}^{
\zeta_\ell(z)-(2m-1)\Theta}
f(z+v',t+v_{2n+1} +2\,{\rm Im}\, (\overline v z))
\,dt
\,dz=\\
&=&
\sum_{
{{\ell\in\Z^{2n},\;\om\cdot\ell=0}\atop{0\leq \ell^{j}
\leq p-1}}
} \int_{H_{\om \ell}}\int_{\zeta_\ell(z) -(2p-1)\Theta}^{
\zeta_\ell(z)+\Theta}
f(z+v',t+v_{2n+1} +2\,{\rm Im}\, (\overline v z))
\,dt
\,dz=\\
&=& 
\sum_{
{{\ell\in\Z^{2n},\;\om\cdot\ell=0
}\atop{0\leq \ell^{j}
\leq p-1}}
} \int_{H_{\om \ell}}\int_{-p\Theta}^{
p\Theta}
f(z+v',t)
\,dt
\,dz=
\\ &=&
\int_{H_\om^p}\int_{-p\Theta}^{p\Theta} f(z+v',t) 
\,dz= \\ &=&
\int_{H_\om^p} g(z+v')\,dz=\\ &=&
\int_{H_\om^p} g(z)\,.
\end{eqnarray*}
By applying the above formula to any~$v$ and to $v=0$
we thus conclude that
$$ \int_{T_v S_\om^p} f=\int_{S_\om^p} f(v\circ \xi)\,d\xi=
\int_{H_\om^p} g(z)= 
\int_{S_\om^p} f\,.$$
This concludes the proof of~(\ref{7uhHhHujJjj6789})
(and so the proof of Lemma~\ref{7ujhHhyhH}).
\end{proof}

We denote by $u_M^p$ the minimal minimizer of
$\Fgot_{S_\om^p}$
under the constraints in ${{\cY}}_M^p$. Obviously, $u_M^1=u_M$.

We now show that changing $S_\om$ into $S_\om^p$
does not change the minimal minimizer:

\begin{lemma}\label{doubling}
$u_M^p=u_M$, for any $p\in\N$, $p\geq 1$.
\end{lemma}

\begin{proof} Let~$u$ be a short hand notation
for~$u_M^p$.
First, we show
that
\begin{equation}\label{yhjk2kk}
u(z,t+2\Theta)=u(z,t)\,,
\end{equation}
for any $(z,t)\in\He^n$.

To this effect, given $a\in\N$,
we define
$$ u_a(\xi)\,:=\, u((0,2a\Theta)\circ\xi)\,.$$
It is easily seen that~$u_a\in {{\cY}}_M^p$.
Also, by Lemma~\ref{7ujhHhyhH}, 
$ \Fgot_{S_\om^p} (u_a)=
\Fgot_{S_\om^p} (u)$, thence~$u_a$
is also a minimizer. Then, as in Lemma~\ref{minmaxiny},
$$ u^\star:=\min \{ u\,,\;u_1\,,\,\dots\,,\; u_{p-1}\}
$$ is a minimizer too.
Since~$u$ is the minimal minimizer, $u\leq u^\star$.
By construction, the converse inequality also holds,
thus $u=u^\star$. Since
$$u_{p-1} ((2\Theta,0)\circ \xi)=u((2p\Theta
,0)\circ \xi)=u(\xi)$$
because~$u\in{{\cY}}_M^p$, we deduce that~$u^\star$
(and, then,~$u$) is~$2\Theta$-periodic in $t$.
This proves~(\ref{yhjk2kk}).

We now show that
\begin{equation}\label{yhjk2kk2}
u((k,0)\circ\xi)=u(\xi)\,,
\end{equation}
for any $k\in\Z^{2n}$ so that $\om\cdot k=0$.
This will be a variation of the proof of~(\ref{yhjk2kk}),
but some care will be needed due to the fact that
the Heisenberg group is non-commutative.
For the proof of~(\ref{yhjk2kk2}),
given $\ell
\in\Z^{2n}$ so that $\om\cdot \ell=0$,
we define
$$ v_\ell(\xi)\,:=\, u((\ell,0)\circ\xi)\,.$$
We have that
\begin{eqnarray}\label{equJJJkL}
v_\ell ((0,2p\Theta)\circ \xi)
=u((\ell,0)\circ(0,2p\Theta)\circ\xi)=
\\
= u((0,2p\Theta)\circ(\ell,0)\circ\xi)=
u(
(\ell,0)\circ\xi)=
v_\ell(\xi)\nonumber
\,,\end{eqnarray}
since $u\in{{\cY}}_M$. Moreover,
if $K\in\Z^{2n}$ is such that $\om\cdot K=0$,
we have that
\begin{eqnarray}\label{eqnar6a6a666ray}&&
\nonumber
(pK,0)\circ (\ell,0) =\\
&=& (0,4p
{\,\rm Im\,}(\overline K \ell ))\circ
(\ell, 0)\circ (pK,0)=\\ &=&
\Big(0,4p \sum_{1\leq a,b\leq 2n-1} K^a \ell^b
{\,\rm Im\,}(k^a k^b )\Big)\circ
(\ell, 0)\circ (pK,0)=\nonumber\\
&=&
\Big(0,4p \sum_{1\leq a,b\leq 2n-1} K^a \ell^b
\Theta_{ab}\Big)\circ
(\ell, 0)\circ (pK,0)\,.\nonumber
\end{eqnarray}
Then, since $\Theta$ divides any $\Theta_{ab}$
(recall the notation 
on page~\pageref{7ujhn8ikj7ujhedfv}),
we gather that
$$
(pK,0)\circ (\ell,0) =
(0,2p\,
\Theta \,\mu)\circ
(\ell, 0)\circ (pK,0)\,,$$
for some $\mu\in \Z$, depending on~$K$ and~$\ell$.

Consequently, using that~$u\in{{\cY}}_M^p$
and~(\ref{yhjk2kk}),
\begin{eqnarray*}
v_\ell((pK,0)\circ \xi)=
u((\ell,0)\circ(pK,0)\circ \xi)=\\
u((pK,0)\circ(\ell,0)\circ \xi)=
u((\ell,0)\circ \xi)=
v_\ell( \xi)
\,.\end{eqnarray*}
This and~(\ref{equJJJkL}) imply that $v_{\ell}\in
{{\cY}}_M^p$.
Moreover,
by Lemma~\ref{7ujhHhyhH}, $ \Fgot_{S_\om^p} (v_{\ell})=
\Fgot_{S_\om^p} (u)$, thence~$v_{\ell}$
is also a minimizer. Then, as pointed out
in Lemma~\ref{minmaxiny},
$$ v^\star:=\min \{ v_{\ell}\,,\;\, 0\leq  \ell^1,\dots,
\ell^{2n-1}\leq p-1\}
$$ is a minimizer too.
Since~$u$ is the minimal minimizer, $u\leq v^\star$.
By construction, the converse inequality also holds,
thus~$u=v^\star$.

We now introduce the following notation.
Given $r\in\Z$, we set $[r]_p\in p\Z$
and $\{r\}_p\in\{1,\dots,p-1\}$ in such a way
$r=[r]_p+\{r\}_p$.

Let now $K\in\Z^{2n}$ be so that $\om\cdot K=0$.
By arguing as in~(\ref{eqnar6a6a666ray}),
we have that
\begin{eqnarray*}
&& (\ell,0) \circ (K,0) =\\
&& =(2b \Theta,0) \circ
\Big( \sum_{j=1}^{2n-1}[(\ell^j+K^j)]_p
\,k^j,0\Big)
\circ
\Big( \sum_{j=1}^{2n-1}\{(\ell^j+K^j)\}_p\,k^j
,0\Big)
\end{eqnarray*}
and so
$$ v_\ell ((K,0)\circ \xi) =
u\Big(
\Big( \sum_{j=1}^{2n-1}\{(\ell^j+K^j)\}_p
\,k^j,0\Big)
\circ\xi
\Big)\,,$$
since $u\in{{\cY}}_M$ and thanks to~(\ref{yhjk2kk2}).

Therefore,
\begin{eqnarray*}
&& v^\star ((K,0)\circ\xi)=\\ &=&
\min \Big\{
u\Big(
\Big( \sum_{j=1}^{2n-1}\{(\ell^j+K^j)\}_p\,
k^j,0\Big)
\circ\xi
\Big)\,,\\ &&
\,\,\;\, 0\leq  \ell^1,\dots,
\ell^{2n-1}\leq p-1\Big\}=
\\&=&
\min \Big\{
u\Big(
\Big( \sum_{j=1}^{2n-1}\tilde\ell^j k^j,0\Big)
\circ\xi
\Big)
\,,\;\, 0\leq  \tilde\ell^1,\dots,
\tilde\ell^{2n-1}\leq p-1\Big\}
=\\
&=& v^\star(\xi)\,.
\end{eqnarray*}
This gives that~$v^\star$ (and thence~$u$)
satisfies~(\ref{yhjk2kk2}),
as desired.

In particular, from~(\ref{yhjk2kk})
and~(\ref{yhjk2kk2}) it follows that
$u=u_M^p\in{{\cY}}_\om$,
and so $\Fgot_{S_\om}(u_M^p)\geq\Fgot_{S_\om}(u_M)$.

On the other hand, it is obvious that $u_M \in {{\cY}}_\om\subseteq
{{\cY}}_\om^p$ and, therefore,
$\Fgot_{S_\om^p}(u_M)\geq\Fgot_{S_\om^p}(u_M)$.

Moreover, if $\nu\in\N$ is the number of
the non-overlapping
components of the union in~(\ref{equasdfgh5tg}) (see also~(\ref{dolla})),
it follows from Lemma~\ref{7ujhHhyhH}
that
$\Fgot_{S_\om^p}=\nu\Fgot_{S_\om}$.

These observations imply that~$u_M$ is a minimizer
of~$\Fgot_{S_\om^p}$
in~${{\cY}}_M^p$ and that~$u_M^p$ is a minimizer
of~$\Fgot_{S_\om}$ in~${{\cY}}_M$.
Since both~$u_M$ and~$u_M^p$
are minimal minimizers,
it follows that $u_M\leq u_M^p\leq u_M$,
as desired.
\end{proof}

Due to Lemma~\ref{doubling}, by possibly replacing~$S_\om$
with~$S_\om^p$,
we now assume that,
given~$r_0$ as in Theorem~\ref{IMPROVE},
the Koranyi balls of radius~$2r_0$
fit inside~$S_\om$. From now on, this suitably large~$p$
will be considered a universal constant, as well as~$r_0$.

In this setting, by a counting argument based on the fact that $\Fgot_{S_\om}$ is bounded independently of $M$, we find
the desired
Koranyi ball
outside the interface:

\begin{lemma}\label{uj8ij88}
Fix $\delta>0$. Fix also $r_0$ suitably large
(according to Theorem~\ref{IMPROVE}).
There exists~$M_0>0$, depending
only on~$\delta$ and on the structural constants of~$\Fgot$,
so that if~$M\geq M_0$, then 
there exists a ball $\BB_\star$ of radius~$r_0$
such that 
$$
\BB_\star \,\subseteq\, S_\om\cap 
\left\{ |(\om,0)\cdot \xi|\leq \frac{M\,|\om|
}{10} \right\}
\cap \{ |u_M| > 1-\delta \}\,.
$$
\end{lemma}

\begin{proof}
Let
$$ \Sigma_M\,:=\,
S_\om\cap 
\left\{ |(\om,0)\cdot \xi|\leq \frac{M
\,|\om|}{10} \right\}\,.$$
We observe that
$$ \Lgot^{2n+1}\left(\Sigma_M
\right)\,\geq \, c_1 \,M\,\Theta\,
|k^1|\,\dots\,|k^{2n-1}|\,,$$
for a suitable $c_1>0$ which may depend
on~$p$ (but it does not depend on either~$\omega$ or~$M$). 
By measure theoretic considerations, we thus deduce
that there are at least N Koranyi balls of radius~$r_0$ 
contained in $\Sigma_M$ whose dialations 
by a constant factor~$\kappa$
(to be chosen appropriately large)
overlap at most a finite number of times (which is independent of~$\om$
and~$M$), with

\begin{equation}\label{8jJJJJKNnN}
N\geq c_2\,M\,\Theta\,|k^1|\,\dots\,|k^{2n-1}|
\end{equation}
The quantity~$c_2$ may depend on~$p$ and~$r_0$,
but not on~$M$ and~$\om$.

Let us suppose that
$\bar n$ of these Koranyi balls
intersect the set $\{ |u_M |< 1-\delta \}$
and let us count how big $\bar n$
can be. For this, if any of the above 
Koranyi balls, say $\BB$, intersects
$\{ |u_M |< 1-\delta \}$ at some point, say $\xi_0$,
then
%% $|u_M (\xi)|\leq 1-\delta/2$ for any $\xi\in \BB_{\e_0}(\xi_0)$,
%% for a small universal~$\e_0>0$, 
%% due to the continuity estimates
%% in~[M95]. Then, 
by Theorem~\ref{IMPROVE},
\begin{equation}\label{7yhgb} 
\Lgot^{2n+1}\Big( \BB_{r_0}(\xi_o) \cap \{
|u_M|\leq 1-\delta\}\Big)\,\geq\, c_3\,,
\end{equation}
for a suitable $c_3$ 
which only depends on $r_0$.

If $\kappa$ is sufficiently large,
$\BB_{r_0}$ is contained in the dialation by 
a factor $\kappa$
of $\BB$, that we denote by $\kappa\BB$.
This and~\eqref{7yhgb} yield that
$$ \Fgot_{\kappa\BB} (u_M) \,\geq \, c_4\,,$$
with $c_4>0$ independent of $\om$
and $M$.
By summing up above all
the
Koranyi balls
intersecting $\{ |u_M |< 1-\delta \}$
(and by using the finite overlapping
properties of these balls),
we thus conclude that
$$ \Fgot_{\Sigma_M}(u_M)\, \geq\, c_5 \bar n\,,$$
with $c_5>0$ independent of~$\om$
and~$M$.

The latter estimate and~\eqref{7ygvftgfFFF}
imply that
$$ \bar n\,\leq\,c_6\,\Theta\,
|k^1|\,\dots\,|k^{2n-1}|\,,$$
where $c_6>0$ is independent of~$\om$
and~$M$.

This and~\eqref{8jJJJJKNnN}
yield that at least one 
Koranyi ball of the family must be outside
$\{ |u_M |< 1-\delta \}$, as long as~$M$ is conveniently
large.
\end{proof}

\subsection{Unconstrained minimizer.}
We are now going to show that, if~$M$ is 
conveniently large, then the constraints 
in~$S_\om$
do not affect the minimal minimizer.

First of all, we show the existence
of a strip of universal size, 
parallel to
the constraints and at a universal distance
from them, on which the minimal minimizer
takes values very close to~$\pm 1$.

\begin{lemma}\label{ujm99ii}
Fix $\delta>0$ and suppose that $M\geq 100$. 
Suppose also that the thesis of
Lemma~\ref{uj8ij88} holds.
Then, there exists \begin{equation}
\label{7uj999}
\lambda \in \left[
-\frac{M}{4}
,\,
\frac{M}{4}
\right]\end{equation}
such that either $u_M(\xi)>1-\delta$
or $u_M(\xi)<-1+\delta$ for any $\xi=(z,t)
\in \He^n$ so that
\begin{equation} \label{ikmh}
z \cdot \frac{\om}{|\om|}\,\in\, [\lambda-1,\,
\lambda+1]
\,.\end{equation}
\end{lemma}

\begin{proof} 
By Lemma~\ref{uj8ij88}, there exists
a  Koranyi
ball of universally large radius
contained in
$\{ |z\cdot w|\leq M/10\}
\cap\{ |u_M |>1-\delta\}$. Then, 
the result 
follows
from the fact that~$u_M$ is continuous, 
$2$-periodic in~$t$ and Birkhoff
(recall Lemmata~\ref{2periodicity} 
and~\ref{uMbir}), 
by applying Lemma~\ref{centrle}.
\end{proof}

We are now in the position
of showing that~$u_M$ does not
vary if we enlarge the constraints,
provided that~$M$ is suitably large:

\begin{prop}\label{ujm99ii-bisco}
Suppose that the hypotheses of Lemma~\ref{ujm99ii}
hold true. Then, $u_{M+a}=u_M$, for any $a\geq 0$.
\end{prop}

\begin{proof} With no loss of
generality, we suppose that $\om_{2n}>0$. 
Let~$H$ be the set of points $\xi=(z,t)\in
\He^n$ \label{7uji880}
satisfying~(\ref{ikmh}) of
Lemma~\ref{ujm99ii}, then
either $u<-1+\delta$ or $u>1-\delta$ in~$H$.

We show, in fact that
\begin{equation}\label{ikkkHHy4567}
{\mbox{$u_M<-1+\delta$ in $H$.}}
\end{equation}
To confirm~(\ref{ikkkHHy4567}), 
we argue by contradiction and
suppose that
$u_M>1-\delta$ in~$H$.
Let $j\in\N$. Then,
by Lemma~\ref{uMbir},
\begin{equation}\label{ujm99} 
T_{(0,\dots,0,j)} 
u_M(\eta ) \geq u_M(\eta) > 1-\delta
\end{equation}
for any $\eta \in H$. 
Take now $\xi=(z,t)$ so that
$\om\cdot z \geq (\lambda +1)|\om|$.
Let $j_0$ be the unique integer so that
$$ j_0\in \left[\frac{|\om|}{\om_{2n}}\Big(
-\lambda -1 +z\cdot \frac{\om}{|\om|}\Big)
,\,
\frac{|\om|}{\om_{2n}}\Big(-
\lambda -1 +z\cdot \frac{\om}{|\om|}\Big)+1
\right)\,.$$
Notice that
$j_0\geq 0$ by construction, thence $j_0\in\N$.
Also, if
$$ \eta:= (0,\dots,0,-j_0,0)\circ \xi\,,
$$
we have that $$\eta=\Big( z-(0,\dots,0,j_0),\,\tau
\Big)\,,
$$ for a suitable $\tau\in\R$ and so
$$ \left( \frac{\om}{|\om|},0
\right)\cdot \eta 
\,\in\, [\lambda-1,\lambda+1]\,,$$
thence $\eta\in H$.

Exploiting~(\ref{ujm99}), we thus conclude
that
$$ u_M(\xi)=
T_{(0,\dots,0,j_0)} 
u_M
(\eta )  >1-\delta\,.
$$
Summarizing, we have just shown that
\begin{equation}
\label{eq4rfg66} u_M
(\xi) > 1-\delta
\end{equation} for any
$\xi=(z,t)\in\He^n$ such that
$\om\cdot z \geq (\lambda +1)|\om|$.

Let now $$ v:= T_{(0,\dots,0,-1)} 
u_M\,.$$
Then, if $(\om,0)\cdot \xi\geq  M|\om|$,
we deduce that
$$ (\om,0)\cdot \Big( 
(0,\dots,0,-1,0)\circ \xi
\Big)=-\om_{2n}+(\om,0)\cdot\xi\geq 
-|\om|+M\,|\om|
\geq (\lambda +1)|\om|\,,$$
due to~(\ref{7uj999}).

Therefore, from~(\ref{eq4rfg66}),
\begin{equation} \label{equa8ijn1}
v(\xi)\geq 1-\delta
\end{equation}
for any $\xi\in\He^n$ so that
$(\om,0)\cdot \xi\geq  M|\om|$.

On the other hand, if
$(\om,0)\cdot \xi\leq  -M|\om|$,
we see that
$$ (\om,0)\cdot \Big( 
(0,\dots,0,-1,0)\circ \xi
\Big)=-\om_n+(\om,0)\cdot\xi
\leq M\,|\om|\,,$$
and so, from the fact that~$u_M\in
\cY_M$, we conclude that
\begin{equation} 
\label{equa8ijn2}
v(\xi)\leq -1+\delta
\end{equation}
for any $\xi\in\He^n$ so that
$(\om,0)\cdot \xi\leq -M|\om|$.

Then,~(\ref{equa8ijn1}), (\ref{equa8ijn2})
Lemma~\ref{yhnUU5tefg}
and Lemma~\ref{2periodicity}
imply that 
\begin{equation} \label{ujm89}
v\in \cY_M\,.\end{equation} 
Moreover, $\Fgot_{S_\om}(u_M)=\Fgot_{S_\om}
(v)$, due to Lemmata~\ref{invfunztra}
and~\ref{yuuiiu}. Consequently, by~(\ref{ujm89}),
$v\in {\mathcal M}_M$ and so, by~(\ref{defmimi}),
$u_M \leq v$.
Since, on the other hand, $u_M \geq v$
by Lemma~\ref{uMbir}, we deduce that $u_M=v$.

Therefore,
$$ u_M (\xi)=u_M \Big( (0,\dots,0,j,0)\circ \xi\Big)\,,$$
for any $j\in \Z$ and any $\xi\in\He^n$.
But then, since $u_M\in \cY_M$,
$$ 1-\delta \leq \lim_{j\rightarrow +\infty}
u_M (0,\dots,0,j,0)=
\lim_{j\rightarrow +\infty}
u_M (0,\dots,0,-j,0)\leq -1+\delta\,,$$
which is a contradiction since we supposed $\delta<1/2$
(see page~\pageref{7yhnjj8}).

This proves~(\ref{ikkkHHy4567}).

\bigskip
We now show that if $\xi=(z,t)\in\He^n$
is so that $\omega\cdot z\leq (\lambda-1)|\om|$
then
\begin{equation}\label{u9999th}
u_M(\xi)<- 1+\delta\,.
\end{equation}
Indeed, if $\xi$ is as above,
let~$j$ be the unique integer so
that
$$ j\in \left(
\frac{|\om|}{\om_{2n}}\Big( 
\lambda-1-
\frac\om{|\om|}
\cdot z\Big)
,\,
\frac{|\om|}{\om_{2n}}\Big(
\lambda-1-
\frac\om{|\om|}
\cdot z\Big)
+1
\right]\,.$$
Then, $j\geq 0$ by construction and so, by
Lemma~\ref{uMbir},
\begin{equation}\label{7uHHHHrRRfF}
T_{(0,\dots,0,j)} u_M \geq u_M\,.
\end{equation}
Moreover, if
\begin{equation}
\label{8ijm} 
\eta\,:= \,
(0,\dots,0,j,0)\circ \xi\,,
\end{equation}
we see that
$$ \left(\frac{\omega}{|\om|},0\right)\cdot \eta\in
[\lambda-1,\,\lambda +1]\,,$$
that is, $\eta\in H$.

Hence,
by~(\ref{ikkkHHy4567}),
$$ T_{(0,\dots,0,j)} u_M(\xi)=
u_M (\eta) <-1+\delta\,.$$
This and~(\ref{7uHHHHrRRfF})
yield the proof of~(\ref{u9999th}).

We now show that 
the lower constraint is irrelevant for~$u_M$.
More precisely, given any $a\geq 0$,
we show that
if
$v$ is $\om$-periodic,
$v(\xi)\leq -1+\delta$ for any $\xi=(z,t)\in\He^n$
so that
$\omega\cdot z\leq -(M+a)|\om|$
and
$v(\xi)\geq 1-\delta$ for any $\xi=(z,t)\in\He^n$
so that
$\omega\cdot z\geq M|\om|$, then
\begin{equation}\label{7ujmk}
\Fgot_{S_\om}(u_M) \,\leq\,
\Fgot_{S_\om}( v)\,.
\end{equation}
To confirm~(\ref{7ujmk}),
take~$v$ as above.
Let $v_a$ be the minimal
minimizer
of $\Fgot_{S_\om}$ in the space
of such $v$'s. Then, by~(\ref{u9999th})
(applied here with the lower constraint
at level~$M+a$
instead of level~$M$), we deduce that
$v_a(\xi) \leq-1+\delta$ for any $\xi\in \He^n$
so that $\omega\cdot z\leq (\lambda-1)|\om|$.
In particular, by~(\ref{7uj999}), we have that $v_a\in
\cY_M$. Therefore, since~$u_M$ is a minimizer
in~$\cY_M$,
$$ \Fgot_{S_\om}(u_M)\leq \Fgot_{S_\om}
(v_a)\leq \Fgot_{S_\om} (v)\,,$$
yielding the proof of~(\ref{7ujmk}).

We are now in the position of completing
the proof of 
Proposition~\ref{ujm99ii-bisco}
by arguing as follows.
Fix $a\geq 0$ and take
\begin{equation}\label{e9ikn1}
w\in \cY_{M+a}\,.
\end{equation}
Take also $j\in\N$ so that 
\begin{equation}\label{e9ikn2}
j> \frac{a\,|\om|}{\om_{2n}}\,.
\end{equation}
Let us define
$$ \tilde w\,:=\, T_{(0,\dots,
0, j)} w\,.$$
Thus, using~(\ref{e9ikn1}) 
and~(\ref{e9ikn2}), we have that
if $\xi\in \He^n$ is so that 
$(\omega,0) \cdot \xi \geq M|\om|$,
then $\tilde w(\xi)\geq 1-\delta$, while
if $(\omega,0)\cdot \xi \leq -(M+j+a)|\om|$, then
$\tilde w(\xi)\leq -1+\delta$.

Therefore, from~(\ref{7ujmk}) 
(used here with~$j+a$
in the place of~$a$), we deduce that
$$ \Fgot_{S_\om} (u_M)\,\leq\,
\Fgot_{S_\om} (\tilde w) \,.$$
Then, by
Lemmata~\ref{invfunztra}
and~\ref{yuuiiu},
$$ \Fgot_{S_\om} (u_M)\,\leq\,
\Fgot_{S_\om} ( w) \,,$$
which completes the proof of
Proposition~\ref{ujm99ii-bisco}.
\end{proof}

We now complete the proof of
Theorem~\ref{classa-per} when $n\geq 2$
(recall~(\ref{n3}))
by arguing as follows.
If $u_\om:=u_M$ is the minimal minimizer constructed above,
we know by~\eqref{consty} and Lemma~\ref{2periodicity}
that the periodicity conditions~\eqref{7uhnHHHB}--\eqref{1.5bis}
and the level set condition~\eqref{7uhnHHH}
hold true.
Also, $u_\om$ satisfies the monotonicity condition
in~\eqref{7uhnHHHBi},
thanks to Lemma~\ref{uMbir} and Definition~\ref{UJNgJhjhg}.

To complete the proof of Theorem~\ref{classa-per} when $n\geq 2$,
we need to show that~$u_\om$ is a local minimizer in any
given ball~$B$.\label{8iujJJjUuHhUY5678}
To this extent, 
we make the following
observation. Let $a\in\Z$ and
$K_1,\dots K_j\in\Z^{2n}$,
with $\om\cdot K_1=\dots=\om\cdot K_j=0$.
Since
$$ (pK_1,0)\circ (pK_j,0)\circ
(0,2ap\Theta)
\in p\Z^{2n}\times p\Theta\Z\,,$$
we deduce that
the distance between the above points
and the origin is large with~$p$ (unless it is zero).
Consequently,
the balls obtained from~$B$ via the group
actions
$$ 
(pK_1,0)\circ (pK_j,0)
\circ(0,2ap\Theta)
$$
are far apart one from the other, for large~$p$.
Therefore, we can extend~$\phi$
in order to fulfill the periodicity
in~(\ref{eqnuuHHy}), that is,
there exists~$\phi^\star\in
{{\cY}}_M^p$
so that~$\phi^\star=\phi$ in~$B$.
Then,
possibly enlarging~$S_\om$ to~$s_\om^p$
as done in Lemma~\ref{doubling}
(note that this does not change the minimal
minimizer), 
we may suppose that~$B$ lies inside~$S_\om$.
Then, using Proposition~\ref{ujm99ii-bisco}, we have that
$\Fgot_{B}(u+\phi)\geq \Fgot_B (u)$, for any $\phi\in C^\infty_0(B)$,
as desired.\label{aqui}

This ends the proof of Theorem~\ref{classa-per} when $n\geq 2$.

We now deal with
the proof of 
Theorem~\ref{classa-per}
when $n=1$.\label{7ujhHHjJkK345}

In this case, the vertical
periodicity of the
minimal minimizer cannot be
recovered from the~$\om$-periodicity,
since the integer base constructed
on page~\pageref{6tgGGYTGG3ed}
boils down to the
two vectors $k^1$ and $k^2=Q\om$ and
so~(\ref{8uhb7}) cannot hold.

To avoid this inconvenience,
the vertical periodicity
needs to be imposed by brute force.
Namely, we substitute
the definition in~(\ref{consty}) with the
following: 
\begin{eqnarray}\label{7ujJJJ}
\nonumber \cY_M &:=&\Big\{u\in S^{1,2}_{\rm loc},
\\ & &\;\,
u\Big( (k^1,0)\circ (z,t)\Big)=
u(z,t)=u(z,t+2)\,\quad\mbox{for any} (z,t)\in \He^1\,,
\nonumber \\
& &\;\, 
u(\xi)\geq 1-\delta\ \mbox{for} \ \xi\cdot (\om,0)\geq
M\,|\om|,\\
& &\;\, u(\xi)\leq -1+\delta\ \mbox{for}\
\xi\cdot (\om,0)\leq -M\,|\om|\Big\}\,.
\nonumber\end{eqnarray}
Then, the arguments on 
pages \pageref{daqui}---\pageref{aqui}
go through verbatim.\footnote{In fact,
in this case, the above arguments
simplify, since $\Theta=1$ and~$u$
is $\om$-periodic if and only if
$$ u\Big( (k^1,0)\circ \xi\Big)=
u(\xi)$$
for any $\xi\in\He^1$. Also,
Lemma~\ref{2periodicity}
is, of course, a trivial
consequence of~(\ref{7ujJJJ}).}
Remarkably, the minimal minimizer obtained
in this way is again a local minimizer
in any domain (that is, the imposed vertical
periodicity does not affect the minimization).
This ends the proof of Theorem~\ref{classa-per}
when~$n=1$.

\begin{rem}
Though we do not make an explicit use of this fact, it is
interesting to note that 
if~$F$ is independent of~$t$, then
so does the minimal minimizer~$u_M$. To confirm this,
fix~$s\in\R$ and
let~$w(z,t):=u_M(z,t+s)$. 
Observe that the map
$$ t\,\longmapsto\, |\gra u_M (z,t)|^2+F\Big(z, u_M (z,t)\Big)$$
is $\Theta$-periodic for any fixed~$z\in \R^{2n}$, due
to~(\ref{chow4}).
Then, the fact that~$F$ does not depend on~$t$
and~(\ref{7uHHHHyHGGG}) imply 
that \begin{eqnarray*} \Fgot_{S_\om}(w)&=&
\int_{H_\om}\int_{-\Theta}^\Theta 
|\gra u_M (z,t+s)|^2+F\Big(z, u_M (z,t+s)\Big)\,dt\,dz\,=\\
&=&
\int_{H_\om}\int_{-\Theta}^\Theta
|\gra u_M (z,t)|^2+F\Big(z, u_M (z,t)\Big)\,dt\,dz\,=\,
\Fgot_{S_\om}(u_M)\,.
\end{eqnarray*}
Thence,
since~$u_M
\in{{\mathcal M}_M}$ by
Lemma~\ref{yhnUUjJ87}, we have that
$w\in{{\mathcal M}_M}$ too. Then, by Lemma~\ref{thJHJ},
$\min\{u_M,w\}
\in{{\mathcal M}_M}$.
So, by~(\ref{defmimi}),
$u_M \leq \min\{u_M,w\}$, which says that $u_M(z,t)\leq u_M (z,t+s)$
for any~$s\in\R$. Then, applying the above observation to
both~$s$ and~$-s$, we get that~$u_M(z,t)=u_M 
(z,t+s)$
for any~$s\in\R$, that is, $u_M$ does not depend on~$t$.

Accordingly, if~$F$ does not depend on~$t$,
the minimal minimizer constructed here agrees with the
uniformly elliptic one of~[V04].
\end{rem}

\section{Proof of Theorem~\ref{classa-qper}}\label{sec:6}

%% First, an easy measure theoretic observation:
%%
%% \begin{lemma}
%% Let $v_j:B\subseteq \R^d\rightarrow\R$ be a sequence
%% of Lebesgue measurable 
%% functions converging to~$v$ uniformly in~$B$.
%% Fix~$c\in\R$. Then,
%% $$ \lim_{j\rightarrow +\infty} \Lgot^d \Big( 
%% B\cap\{ v_j =c\neq v\}
%% \Big)\,=\,0\,.$$
%% \end{lemma}
%%
%% \begin{proof} Possibly replacing~$v_j$
%% by~$v_j-c$, we may and do assume that~$c=0$.
%% Let $\e_j>0$ be an infinitesimal decreasing sequence
%% so that
%% $$ \sup_B |v_j-v|\leq \e_j\,.$$
%% Then,
%% $$ \{ v_j =0\neq v\}\subseteq \{ 0<|v|\leq \e_j\}
%% =: E_j\,.$$
%% By construction, $E_{ j+1}\subseteq E_j$
%% and
%% $$\bigcap_{j\in\N } E_j\,=\,\emptyset\,.$$
%% Then, by the Monotone Convergence
%% Theorem,
%% $$ \lim_{j\rightarrow +\infty} \Lgot^d ( 
%% E_j)\,=\,0\,,$$
%% which yields the desired claim.
%% \end{proof}

We 
consider~$\om\in\R^{2n}\setminus\Q^{2n}$ 
and~$\om_j\in\Q^{2n}$
be a sequence of rational vectors
approaching~$\om$.
Exploiting Theorem~\ref{classa-per}, we obtain the existence of a 
function~$u_j:\He^n\rightarrow[-1,1]$ 
which is a
local minimizer for~$\Fgot$ in any bounded
domain, is $\om_j$-periodic and has the level
sets $\{|u_j|\leq 1-\delta\}$ trapped inside a slab
normal to $(\om_j,0)$
of size~$M_0$. We stress that~$M_0$
may depend on~$\delta$ and on the structural
constants of the problem, but it is independent
of~$j$. Then, to complete the
proof of Theorem~\ref{classa-qper},
we need to show that, perhaps passing 
to a 
subsequence,~$u_j$ 
converges locally uniformly
to a suitable~$u$
which is a local minimizer in any bounded domain.

For this, we observe that~$|u_j|\leq 1$ by construction
and that the modulus of continuity of~$u_j$
is uniformly bounded, thanks to [M95].
Therefore, up to subsequences,
$u_j$ converges 
locally uniformly
to a suitable~$u$.

Given a bounded domain~$\Omega$,
we now show that~$u$ is a local minimizer
of~$\Fgot_\Omega$. For this, let~$\phi\in C^\infty_0 (\Omega)$
and~$R>0$ so large that~$\Omega \subset \BB_R$.

Observe that, possibly taking disjoint unions of~$S_\om$
as in Lemma~\ref{doubling} (which does not affect the minimal
minimizer), we may assume that~$\BB_{R+2}$ sits inside~$S_\om$.

Then,~$u_j$ is a local minimizer in~$\BB_{R+2}$.
Accordingly, there exists a suitable~$C_R>0$ 
so that 
\begin{equation}\label{lduest} 
\int_{\BB_{R+1}} |\gra u_j|^2\,\leq\,
C_R\end{equation}
for any $j\in\N$, due to the Caccioppoli-type
estimate in~(\ref{cacci}).
Thus, possibly taking subsequences,
we have that~$\gra u_j$ converges to~$\gra u$
weakly in~$L^2(\BB_{R+1})$, due to the embedding results
in~[VSCC92].

Therefore,
\begin{equation}\label{kin eq}
\int_{\BB_{R+\alpha}
}|\gra u|^2\,\leq\,
\liminf_{j\rightarrow+\infty} \int_{
\BB_{R+\alpha}
}|\gra u_j|^2\,,
\end{equation}
for any~$\alpha\in[0,1]$.

Also,
by the continuity of~$F$
when the second variable is in~$(-1,1)$,
we deduce from the
uniform convergence
of~$u_j$ and Fatou's Lemma
that
\begin{equation}\label{ftu eq}
\begin{split}
& \liminf_{j\rightarrow+\infty}
\int_{\BB_{R}} F(\xi, u_j(\xi))\,d\xi \geq\\
\geq &
\liminf_{j\rightarrow+\infty}
\int_{\BB_{R}\cap
\{ |u|< 1\}
} F(\xi, u_j(\xi))\,d\xi \geq\\
\geq &
\int_{\BB_{R}\cap\{ |u|< 1\}
} F(\xi, u(\xi))\,d\xi=\\
=&
\int_{\BB_{R}
} F(\xi, u (\xi))\,d\xi\,.
\end{split}
\end{equation}
By collecting the estimates in~(\ref{kin eq})
and~(\ref{ftu eq}), we thus conclude that
\begin{equation}\label{7uhn8uj}
\Fgot_{\BB_{R}} (u)\,\leq\,\liminf_{j\rightarrow
+\infty}\Fgot_{\BB_{R}} (u_j)\,.
\end{equation}
We now 
define
\begin{equation} \label{eqyhj001}
\e_j \,:=\,\sup_{\BB_{R+1}}
|u-u_j|\,.\end{equation}
By the uniform convergence
of~$u_j$, we have that~$\e_j$ converges
to zero. 
Moreover, by~(\ref{lduest})
and~(\ref{kin eq}),
\begin{equation}
\label{3.bis} 
\int_{\BB_{R+1}}|\gra u|^2 \leq C_R
\end{equation}
and so, by the absolute continuity of
the Lebesgue integral,
\begin{equation}\label{tyhj7u}
\lim_{j\rightarrow +\infty}
\int_{\BB_{R+\e_j}\setminus
\BB_R}|\gra u|^2\,=\,0\,.
\end{equation}
We also set
\begin{equation}\label{muggei}
\mu_j\,:=\, \Lgot^{2n+1}(\BB_{R+\e_j}
\setminus\BB_R)\,,
\end{equation}
thence~$\mu_j$ also
converges
to zero. 

Let now~$\tilde\tau_j\in 
C^\infty_0([0,R+\e_j])$
so that $0\leq \tilde\tau_j\leq 1$, $\tau_j(t)=1$
for any $t \in [0,R]$ and $|\tilde\tau_j'|\leq 10/\e_j$.
Let~$\tau_j$ be the radial (with respect to the
Koranyi ball) function defined by~$\tau_j
(\xi):=\tilde\tau_j(\rho)$
(recall the notation in~(\ref{radialvars})).
Then, by~(\ref{grast}), 
\begin{equation}\label{eqyhj002}
|\gra \tau_j|\leq 
\frac{{{\,\rm const\,}}}{\e_j}\,.\end{equation}
Let also
$$ w_j\,:=\, \tau_j u + (1-\tau_j) u_j+\phi\,.$$
By construction,
$w_j(\xi) =u_j(\xi)$ for any $\xi \in \He^n\setminus\BB_{
R+\e_j}$, and so, by the minimizing property of~$u_j$,
\begin{equation}\label{ujkooo998}
\begin{split}
&\Fgot_{\BB_R} (u_j) +\int_{\BB_{R+\e_j}
\setminus \BB_R }|\gra u_j|^2\leq 
\\ \leq&
\Fgot_{\BB_{R+\e_j}}(u_j) \leq \\ 
\leq& \Fgot_{\BB_{R+\e_j}}
(w_j)
= \\=&
\Fgot_{\BB_{R}}
(w_j)
+\Fgot_{\BB_{R+\e_j}\setminus \BB_R}
(w_j)=\\
=&
\Fgot_{\BB_{R}}
(u+\phi)
+\Fgot_{\BB_{R+\e_j}\setminus\BB_R}
(w_j)
\,.\end{split}
\end{equation}
On the other hand,
exploiting~(\ref{eqyhj002}) and~(\ref{eqyhj001}),

\begin{eqnarray*}
|\gra w_j| &=&
|\gra\tau_j (u-u_j) +\tau_j\gra u+(1-\tau_j)
\gra u_j+\gra \phi|\leq\\
&\leq& \frac{C}{\e_j}|u-u_j| +|\gra u|+|\gra u_j|+|\gra \phi|\\
&\leq& C+|\gra u|+|\gra u_j|
\end{eqnarray*}
for a suitable~$C>0$ which may depend
on~$R$ and~$\phi$, but it is
independent of~$j$.

Using the above estimate, Cauchy Inequality 
and~(\ref{muggei}),
we obtain that
\begin{equation}\label{ujnkkujnasdf}
\begin{split}
&\Fgot_{\BB_{R+\e_j}\setminus\BB_R}
(w_j) \leq \\ &\leq \int_{\BB_{R+\e_j}\setminus\BB_R}
\Big(|\gra u|^2+|\gra u_j|^2\Big)+\\
&+C \sqrt{\mu_j}\,\left( 
\sqrt{\int_{\BB_{R+1}} |\gra u|^2}+
\sqrt{\int_{\BB_{R+1}} |\gra u_j|^2}+1
\right)+\\ &+C
\sqrt{\int_{\BB_{R+1}} |\gra u_j|^2\,
\int_{\BB_{R+\e_j}\setminus\BB_R} |\gra u|^2}
\,.\end{split}
\end{equation}
We deduce from~(\ref{lduest}),~(\ref{3.bis}),~(\ref{tyhj7u})
and~(\ref{ujnkkujnasdf}) that
$$ \Fgot_{\BB_{R+\e_j}\setminus\BB_R}
(w_j) \leq \delta_j+
\int_{\BB_{R+\e_j}\setminus\BB_R}
|\gra u_j|^2\,,$$
for a suitable $\delta_j$ converging to zero.

Consequently, from~(\ref{ujkooo998}),
\begin{equation} \label{6yhuu88745t}
\Fgot_{\BB_R} (u_j) \leq 
\Fgot_{\BB_R}(u+\phi) +\delta_j\,.
\end{equation}
By~(\ref{7uhn8uj}),~(\ref{6yhuu88745t})
and the fact that~$\delta_j$ converges to zero,
we thus conclude that
$$ \Fgot_{\BB_R} (u)\,\leq\,
\Fgot_{\BB_R} (u+\phi)\,.$$
Thence, since~$\phi$ vanishes outside~$\Omega$,
$$ \Fgot_{\Omega} (u)\,\leq\,
\Fgot_{\Omega} (u+\phi)\,,$$
which shows that~$u$ is a local minimizer
for~$\Fgot$ in~$\Omega$.
This completes the proof of
Theorem~\ref{classa-qper}.

\section{Proof of Theorem~\ref{backto}}\label{sec:7}

We suppose that $\om\in\Q^{2n}\setminus\{0\}$
(the general case then follows by a standard
limiting
argument, see [CL01]).

Let
$Q(\xi):=\alpha(\xi)^2$
and, for $N\in\N$, 
$$ \Fgot^{(N)}_\Omega (u)\,:=\,
\int_\Omega \frac 1 N\,
|\gra u(\xi)|^2 \,+\,
N\,Q(\xi) (1-u(\xi)^2)^2\,
d\xi\,.$$
Let~$u_N$ be
the minimal minimizer of~$\Fgot^{(N)}$,
as defined in~\eqref{defmimi},
under
the constraint in~\eqref{consty}.
In particular,
\begin{equation}\label{uMbuu4}
\{ |u|< 1-\delta\}\,\subseteq\,\{
|(\om,0)\cdot\xi|\leq M
\}
\,.\end{equation}
Then, by Lemma~\ref{uMbir},
\begin{equation}\label{uMbuu5}
{\mbox{
$u_N$ satisfies the Birkhoff property.}}
\end{equation}
We also set
$$ v_N(z,t)\,:=\, u_N\left(\frac z N, \frac{t}{N^2}
\right)\,.$$
Then, $v_N$ minimizes
\begin{equation}\label{yhbnHJHHHY}
\int |\gra v(z,t)|^2 +Q \left(\frac z N, \frac{t}{N^2}
\right)\, \Big(1-v(z,t)^2\Big)^2\,dz\,dt
\end{equation}
under
the constraint
that $v_N (\xi)\geq 1-\delta$ if 
$\xi\cdot (\om,0)\geq MN$ and
$v_N(\xi)\leq -1+\delta$ if
$\xi\cdot (\om,0)\leq -MN$. Also,
if~$k\in\Z^{2n}$ is so that $\om\cdot k=0$,
then $v_N ((Nk,0)\circ\xi)=v_N(\xi)$ for
any~$\xi\in\He^n$.

Notice that the structural constants of the functional 
in~\eqref{yhbnHJHHHY} (as they appear on page~\pageref{7HHHGGHHhj11})
may be bounded independently of~$N$.

Then, by performing the counting argument
in Lemma~\ref{uj8ij88}, we find a Koranyi
ball of radius~$N r_0$
contained in
$$\{|\xi\cdot (\om,0)|\leq MN\}\cap
\{ |v_N|>1-\delta\}\,,$$
if~$M\geq M_0 |\om|$ and~$M_0$
is large enough.

By scaling
back, we 
find a Koranyi
ball of radius~$r_0$
contained in
$$\{|\xi\cdot (\om,0)|\leq M\}\cap
\{ |u_N|>1-\delta\}\,.$$
Then, by using~\eqref{uMbuu5}
and the argument
in Proposition~\ref{ujm99ii-bisco},
we deduce that~$u_N$ is a free
minimizer.

Then, possibly taking subsequences, using the results in~[MSC01], 
we have 
that~$u_N$ 
converges in~$L^1_{\rm loc}$ and almost everywhere to a step
function $\chi_{E}-\chi_{\He^n \setminus E}$
and $\partial E$ has minimal surface area
with respect to the metric~$\alpha$.
This and~\eqref{uMbuu4}
imply that $$\partial E\,\subseteq
\,\{|\xi\cdot (\om,0)|\leq M\}
$$ 
as long as $M\geq M_0 |\om|$ and $M_0$
is universally large, ending the proof of
Theorem~\ref{backto}.

\bigskip\bigskip

{\footnotesize

\section*{References}

\begin{itemize}

\item[{[AAC01]}]
Alberti, G.; Ambrosio, L.; Cabr\'e, X. -- {\em
On a long-standing
conjecture of E. De~Giorgi: symmetry in 3D for general nonlinearities and
a local minimality property},
Acta Appl. Math.  65
(2001),  no. 1-3, 9--33.

\item[{[AC81]}] Alt, H. W.; Caffarelli, L. A. -- {\em
Existence and regularity for a minimum problem with free boundary},
J. Reine Angew. Math. 325 (1981),
105--144.

\item[{[ACF84]}] Alt, H. W.; Caffarelli, L. A.;
Friedman, A. -- {\em
A free boundary problem for quasilinear
elliptic equations},
Ann. Scuola Norm. Sup. Pisa Cl. Sci. (4) 11 (1984), no. 1, 1--44.

\item[{[B90]}]
Bangert, V. -- {\em Laminations of {$3$}-tori by least area 
surfaces}, Analysis,
et cetera, Academic Press, Boston, MA, 1990, pp.~85--114.

\item[{[B03]}]
Birindelli, I. -- {\em
Superharmonic functions in the Heisenberg group: estimates and 
Liouville theorems},
NoDEA, Nonlinear Differ. Equ. Appl. 10 (2003), no. 2, 171--185.

\item[{[BL03]}]
Birindelli, I.; Lanconelli, E. -- {\em
A negative answer to a one-dimensional symmetry problem in the Heisenberg 
group},
Calc. Var. Partial Differ. Equ. 18 (2003), no. 4, 357--372.

\item[{[CC95]}] Caffarelli, L. A.; C\'ordoba, A. -- {\em
Uniform convergence of a singular perturbation problem},
Comm. Pure Appl. Math. 48 (1995), no.
1, 1--12. 

\item[{[CL01]}]
Caffarelli, L. A.; de~la Llave, R. -- {\em Planelike minimizers in
periodic media}, Comm. Pure Appl. Math. 54 (2001), no.~12,
1403--1441.

\item[{[CL05]}]
Caffarelli, L. A.; de~la Llave, R. --
{\em Interfaces of ground 
states in
{I}sing models with periodic coefficients}, J. Stat. Phys. {118}
(2005), no.~3-4, 687--719.

\item[{[CL98]}]
Candel A.; de~la Llave, R. -- {\em On the {A}ubry-{M}ather theory in
statistical mechanics}, Comm. Math. Phys. 192 (1998), no.~3,
649--669.

\item[{[DG79]}] De~Giorgi, E. -- {\em
Convergence problems for functionals and operators},
Proceedings of the International Meeting on Recent Methods in Nonlinear
Analysis (Rome, 1978), 131--188,
Pitagora, Bologna, 1979.

\item[{[DN05]}] 
Dynnikov, I.A.; Novikov, S.P. -- {\em
Topology of quasiperiodic functions on the plane},
Russ. Math. Surv. 60 (2005), no.~1, 1--26.

\item[{[FS74]}] Folland, G. B.; Stein, E. M. -- {\em 
Estimates for the 
$\bar \partial \sb{b}$ 
complex and analysis on the Heisenberg group},
Comm. Pure Appl. Math. 27  
(1974), 429--522. 

\item[{[FSSC01]}]
Franchi, B.; Serapioni,~R.; Serra Cassano,~F.~-- {\em
Rectifiability and perimeter in the Heisenberg group},
Math. Ann. 321 (2001), no.~3, 479--531.

\item[{[GL58]}] Ginzburg, V.; Landau, L. -- {\em
On the theory of superconductivity}, Zh. \`Eksper. Teoret. Fiz. 20
(1950), 1064--1082.

\item[{[GP58]}] Ginzburg, V. L.; Pitaevski\u\i, L. P. -- {\em
On the theory of superfluidity}, Soviet
Physics. JETP 34/7 (1958), 858--861.

\item[{[G84]}]
Giusti, E. -- {\em 
Minimal surfaces and functions of bounded variation}, 
Monographs in Mathematics, 80. 
Birkh\"auser Verlag, Basel-Boston, Mass., 1984.

\item[{[G03]}]
Giusti, E. -- {\em
Direct methods in the calculus of variations},
World Scientific Publishing 
Co., Inc., River Edge, NJ, 2003. viii+403 pp.

\item[{[G85]}]
Gurtin, M. E. -- {\em
On a theory of phase transitions with interfacial energy}, 
Arch. Rational Mech. Anal. 87 (1985), no. 3, 187--212. 

\item[{[H32]}]
Hedlund, G. A. -- {\em Geodesics on a two-dimensional {R}iemannian 
manifold
with periodic coefficients}, Ann. of Math. 33 (1932), 719--739.

\item[{[M95]}]
Marchi, S. -- {\em H\"older continuity 
and Harnack inequality for De Giorgi 
classes related to H\"ormander vector fields},
Ann. Mat. Pura Appl.
168  
(1995), no. 10, 171--188. 

\item[{[M87]}]
Modica, L. -- {\em
The gradient theory of phase transitions
and the minimal interface criterion},
Arch. Rational Mech. Anal. 98 (1987), no.
2, 123--142. 

\item[{[MSC01]}]
Monti, R.; Serra Cassano, F. -- {\em
Surface measures in Carnot-Carathéodory spaces},
Calc. Var. Partial Differ. Equ. 13 (2001), no. 3, 339--376.

\item[{[M86]}]
Moser, J. -- {\em Minimal solutions of variational problems on a 
torus},
Ann. Inst. H. Poincar\'e Anal. Non Lin\'eaire 3 (1986), no.~3,
229--272.

\item[{[R79]}] Rowlinson, J. S. -- {\em
Translation of J. D. van der Waals' ``The thermodynamic theory of
capillarity under the hypothesis of a continuous variation of density''},
J. Statist. Phys. 20 (1979), no. 2,
197--244. 

\item[{[PV05a]}] Petrosyan, A.; Valdinoci, E. --
{\em Geometric properties of
Bernoulli-type minimizers},
Interfaces Free Bound. 7 (2005), 55--78.

\item[{[PV05b]}] Petrosyan, A.; Valdinoci, E. --
{\em Density estimates for a degenerate/singular
Ginz\-burg-Landau model},
SIAM J. Math. Anal. 36 (2005), no. 4, 1057--1079.

\item[{[RS03]}] Rabinowitz, P. H.; Stredulinsky, E. --
{\em Mixed states 
for an
{A}llen-{C}ahn type equation}, Comm. Pure Appl. Math. 56 
(2003),
no.~8, 1078--1134.

\item[{[RS04]}] Rabinowitz, P. H.; Stredulinsky, E. --
{\em On some results of {M}oser and 
of
{B}angert}, Ann. Inst. H. Poincar\'e Anal. Non Lin\'eaire 21 
(2004),
no.~5, 673--688.

\item[{[V04]}] Valdinoci, E. --
{\em Plane-like minimizers in periodic media: jet flows and
Ginzburg-Landau models},
J. Reine Angew. Math. 
574 (2004), 147--185.

\item[{[VSCC92]}] Varopoulos,  N. Th.;
Saloff-Coste, L.; Coulhon, T. --
{\em
Analysis and geometry on groups  
},
Cambridge University Press, UK, 1992. xii+156 pp.
\end{itemize}

\bigskip\bigskip
$$ \star \qquad\star\qquad\star$$
\bigskip\bigskip

\end{document}